\documentclass[12pt]{amsart}
\usepackage{amsmath, amssymb}

\theoremstyle{plain}
\newtheorem{theorem}{Theorem}[section]
\newtheorem{proposition}[theorem]{Proposition}
\newtheorem{corollary}[theorem]{Corollary}
\newtheorem{lemma}[theorem]{Lemma}

\theoremstyle{definition}
\newtheorem{definition}[theorem]{Definition}

\theoremstyle{remark}
\newtheorem*{remark}{Remark}

\numberwithin{equation}{section}

\def \remarks {\noindent {\bf Remarks.}\ \ }

\def \N {\mathbb{N}}
\def \R {\mathbb{R}}

\def \Z {\mathbb{Z}}
\def \E {\mathbb{E}}
\def \P {\mathbb{P}}

\def \one {{\bf 1}}
\def \DD {\mathcal{D}}

\def \NN {\mathcal{N}}

\def \PP {\mathcal{P}}

\def \a {\alpha}
\def \b {\beta}

\def \e {\varepsilon}

\def \d {\delta}

\def \k {\kappa}
\def \l {\lambda}

\def \s {\sigma}

\def \z {\zeta}

\def \< {\langle}
\def \> {\rangle}
\def \^ {\widehat}

\def \dist {{\rm dist}}

\def \Span {{\rm span}}

\def \dens {{\rm dens}}
\def \supp {{\rm supp}}
\def \Sparse {{\mathit{Sparse}}}
\def \Comp {{\mathit{Comp}}}
\def \Incomp {{\mathit{Incomp}}}

\newcommand{\norm}[1]{\left \| #1 \right \|}
\def \etc {,\ldots,}

\begin{document}
\title[]{The Littlewood-Offord Problem and invertibility of random matrices}

\author{Mark Rudelson
  \and Roman Vershynin}

\thanks{\hspace{-.5cm} M.R.: Department of Mathematics,
  University of Missouri,
  Columbia, MO 65211, USA. \\
  e-mail: rudelson@math.missouri.edu   \\
  R.V.: Department of Mathematics,
  University of California,
  Davis, CA 95616, USA.\\
  e-mail: vershynin@math.ucdavis.edu  \\
  \newline
  M.R. was supported by NSF DMS grant 0556151.\\
  R.V. was supported by the Alfred P.~Sloan Foundation
  and by NSF DMS grants 0401032 and 0652617.}




\maketitle

\begin{abstract}
  We prove two basic conjectures on the distribution of the smallest
  singular value of random $n \times n$ matrices with independent entries.
  Under minimal moment assumptions, we show that the
  smallest singular value is of order $n^{-1/2}$, which is optimal for Gaussian
  matrices. Moreover, we give a optimal estimate on the tail probability.
  This comes as a consequence of a new
  and essentially sharp estimate in the Littlewood-Offord problem:
  for i.i.d. random variables $X_k$ and real numbers $a_k$,
  determine the probability $p$ that the sum $\sum_k a_k X_k$ lies
  near some number $v$.
  For arbitrary coefficients $a_k$ of the same order of magnitude,
  we show that they essentially lie in an arithmetic progression of length
  $1/p$.
\end{abstract}

\section{Introduction}

\subsection{Invertibility of random matrices}

In this paper we solve two open problems on the distribution of the
smallest singular value of random matrices.

Let $A$ be an $n \times n$ matrix with real or complex entries. The
{\em singular values} $s_k(A)$ of $A$ are the eigenvalues of $|A| =
\sqrt{A^*A}$ arranged in the non-increasing order. Of particular
significance are the largest and the smallest singular values
$$
s_1(A) = \sup_{x:\; \|x\|_2 = 1} \|Ax\|_2, \qquad s_n(A) =
\inf_{x:\; \|x\|_2 = 1} \|Ax\|_2.
$$
These quantities can obviously be expressed in terms of the spectral
norm -- the operator norm of $A$ considered as an operator on $\ell_2^n$.
Indeed, $s_1(A) = \|A\|$, and if the matrix $A$ is non-singular then
$s_n(A)=1/\norm{A^{-1}}$. The smallest singular value thus equals
the distance from $A$ to the set of singular matrices in the
spectral norm.

The behavior of the largest singular value of random matrices $A$
with i.i.d. entries is well studied.
The weakest assumption for its regular behavior is
boundedness of the {\em fourth moment} of the entries; then
\begin{equation}                \label{lmax asymptotic}
  s_1(A) \sim n^{1/2} \quad \text{with high probability}.
\end{equation}
Indeed, by \cite{YBK, BSY} the finite fourth moment is {\em necessary
and sufficient} for $s_1(A)/n^{1/2}$ to have an almost sure
limit as $n \to \infty$, and this limit equals $2$.
Latala \cite{Lat} showed that \eqref{lmax asymptotic} holds
under the forth moment assumption even if entries are not
identically distributed.

Much less has been known about the behavior of the smallest singular
value. In the classic work on numerical inversion of large matrices,
von~Neumann and his associates used random matrices to test their
algorithms, and they speculated that
\begin{equation}                \label{whp}
  s_n(A) \sim n^{-1/2} \quad \text{with high probability}
\end{equation}
(see \cite{vN}, pp.~14, 477, 555). In a more precise form, this
estimate was conjectured by Smale \cite{S 85} and proved by Edelman
\cite{E 88}  and Szarek \cite{Sz} for {\em random Gaussian matrices}
$A$, those with i.i.d. standard normal entries. Edelman's theorem
states that for every $\e \ge 0$
\begin{equation}                \label{tail}
  \P \big( s_n(A) \le \e n^{-1/2} \big) \sim \e.
\end{equation}

Prediction \eqref{whp} for general random matrices has been an open
problem, unknown even for the {\em random sign matrices} $A$,
those whose entries are $\pm 1$ symmetric random variables.
In this paper we prove the prediction \eqref{whp} in full generality
under the aforementioned fourth moment assumption.

\begin{theorem}[Invertibility: fourth moment]       \label{t: 4}
  Let $A$ be an $n \times n$ matrix whose entries are independent
  real random variables with variances at least $1$ and fourth moments
  bounded by $B$.
  Then, for every $\d > 0$ there exist $\e > 0$ and $n_0$
  which depend (polynomially) only on $\d$ and $B$, and such that
  $$
  \P \big( s_n(A) \le \e n^{-1/2} \big) \le \d
    \qquad \text{for all $n \ge n_0$}.
  $$
\end{theorem}

This shows in particular that the median of $s_n(A)$ is of order
$n^{-1/2}$.

\medskip

Under stronger moment assumptions, more is known about the distribution
of the largest singular value, and similarly one hopes to know more
about the smallest singular value.

Indeed, Soshnikov \cite{So} proved that the limiting distribution of $s_1(A)$
is precisely the Tracy-Widom law for all matrices with i.i.d. subgaussian entries.
Recall that a random variable $\xi$ is called {\em subgaussian} if its tail is
dominated by that of the standard normal random variable:
there exists $B > 0$ such that
\begin{equation}   \label{subgaussian}
  \P (|\xi| > t) \le 2 \exp(-t^2/B^2) \qquad \text{for all $t > 0$}.
\end{equation}
The minimal $B$ here is called the {\em subgaussian moment}\footnote{
  In the literature in geometric functional analysis, the subgaussian
  moment is often called the $\psi_2$-norm.}
of $\xi$.
Inequality \eqref{subgaussian} is
often equivalently formulated as a moment condition
\begin{equation}    \label{subgaussian moments}
   (\E |\xi|^p )^{1/p} \le C B \sqrt{p} \qquad \text{for all $p \ge 1$},
\end{equation}
where $C$ is an absolute constant.
The class of subgaussian random variables includes many random variables
that arise naturally in applications, such as normal, symmetric $\pm 1$,
and in general all bounded random variables.

One might then expect that the estimate \eqref{tail} for the
distribution of the smallest singular value of Gaussian matrices
should hold for all subgaussian matrices. Note however that
\eqref{tail} fails for the random sign-matrices, since they are
singular with positive probability. Estimating the singularity
probability for random sign-matrices is a longstanding open problem.
Even proving that it converges to $0$ as $n \to \infty$ is a
nontrivial result due to Koml\'os \cite{K 67}. Later Kahn, Koml\'os
and Szemer\'edi \cite{KKS} showed that it is exponentially small:
\begin{equation}            \label{KKS}
  \P \big( \text{random sign matrix $A$ is singular} \big) < c^n
\end{equation}
for some universal constant $c \in (0,1)$. The often conjectured
optimal value of $c$ is $1/2 + o(1)$ \cite{O 88, KKS}, and
the best known value $3/4 + o(1)$ is due to Tao
and Vu \cite{TV det, TV singularity}.

Spielman and Teng \cite{ST} conjectured that \eqref{tail} should hold
for the random sign matrices up to an exponentially
small term that accounts for their singularity probability:
$$
\P \big( s_n(A) \le \e n^{-1/2} \big) \le \e + c^n.
$$

In this paper, we prove Spielman-Teng's conjecture for all matrices
with subgaussian i.i.d. entries, and up to a constant factor
which depends only on the subgaussian moment.

\begin{theorem}[Invertibility: subgaussian]   \label{t: subgaussian}
  Let $\xi_1 \etc \xi_n$ be independent centered real random variables
  with variances at least $1$ and subgaussian moments bounded by $B$.
  Let $A$ be an $n \times n$ matrix whose
  rows are independent copies of the random vector $(\xi_1,\ldots,\xi_n)$.
  Then for every $\e \ge 0$ one has
  \begin{equation}                      \label{eq main}
    \P \big( s_n(A) \le \e n^{-1/2} \big)
    \le C\e +  c^n,
  \end{equation}
  where $C > 0$  and $c \in (0,1)$ depend (polynomially) only on $B$.
\end{theorem}

\begin{remarks}
1. For $\e = 0$, Theorem~\ref{t: subgaussian} yields an exponential bound
for the singularity probability:
$$
\P \big( \text{random subgaussian matrix $A$ is singular} \big) < c^n.
$$
Thus Kahn-Koml\'os-Szemer\'edi's bound \eqref{KKS} holds
for all subgaussian matrices.
Moreover, while \eqref{KKS} estimates the probability that a random
matrix {\em belongs} to the set of singular matrices,
Theorem~\ref{t: subgaussian} estimates the {\em distance} to that set.

2. The bounds in  Theorem~\ref{t: subgaussian} are precise.
Edelman's bound \eqref{tail} shows that the term $\e n^{-1/2}$ is
optimal for the Gaussian matrix, while the term $c^n$ is optimal for
a random sign-matrix.

3. For simplicity, we state and prove all our results over the real field. 
However, our arguments easily generalize to the complex field; 
see e.g. \cite{PZ}.

 4. A weaker result was recently proved by the first author
\cite{R} who showed that $\P \big( s_n(A) \le \e n^{-3/2} \big) \le
C\e + C n^{-1/2}$. He later  improved the term $n^{-1/2}$ to $c^n$.
Shortly after that, both authors of this paper independently
discovered how to reduce the term $n^{-3/2}$ to the sharp order
$n^{-1/2}$. In December 2006, the second author found a new way to
prove the sharp invertibility estimate by obtaining an essentially
optimal result for the Littlewood-Offord problem as stated in
Theorem~\ref{t: small ball intro}. We thus decided to publish
jointly, and some of the arguments were improved during the final
stage of our work.

5. Another weaker result was recently proved by Tao and Vu \cite{TV}
for random sign matrices. They showed that for every $A>0$ there
exists $B>0$ such that $s_n(A) \ge n^{-B}$ holds with probability $1
- O_A(n^{-A})$.
\end{remarks}

\subsection{The Littlewood-Offord problem}

Our results on random matrices come as a consequence of a new and
essentially sharp estimate in the Littlewood-Offord problem
\cite{B}, \cite{FF}. A classical theme in Probability Theory is the
study of the random sums
\begin{equation}                    \label{random sum}
  S := \sum_{k=1}^n a_k \xi_k,
\end{equation}
where $\xi_1,\ldots,\xi_n$ are independent identically distributed
random variables and $a = (a_1,\ldots,a_n)$ is a vector of real
coefficients.

The large deviation theory demonstrates that $S$ nicely concentrates
around its mean. On the other hand, by the central limit theorem,
one can not expect tighter concentration than that of the
appropriately scaled Gaussian random variable. However, rigorous
anti-concentration estimates are hard to prove (see \cite{LS}),
especially for discrete random variables $\xi_k$. The
Littlewood-Offord problem  thus asks to estimate the {\em small ball
probability}
$$
p_\e(a) := \sup_{v \in \R} \; \P (|S-v| \le \e).
$$
A small value of $p_\e(a)$ would mean that the random sums $S$ are
well spread.

For the {\em random Gaussian sums}, i.e. for $\xi_k$ being standard
normal random variables, the small ball probability for each $\e$
depends only on the Euclidean norm of the coefficient vector $a$ and
not on its direction, and one has $p_\e(a) \sim \e/\|a\|_2$.

For most other distributions, $p_\e(a)$ depends on the direction of
$a$, and determining the asymptotics is hard. A remarkable and
extensively studied case is for the {\em random sign-sums} $\sum \pm
a_k$, i.e. for symmetric $\pm 1$ random variables $\xi_k$. The small
ball probability strongly depends on the direction of the
coefficient vector: for example, $p_0(a) = 1/2$ for $a =
(1,1,0,\ldots,0)$ while $p_0(a) \sim n^{-1/2}$ for $a =
(1,1,\ldots,1)$.

The coefficient vectors with few non-zero coordinates turn out to be
the only obstacle for nontrivial estimates on the small ball
probability. The classical result of Littlewood and Offord
strengthened by Erd\"os \cite{E 45} states that if all $|a_k| \ge 1$
then for the random sign-sums one has
\begin{equation}                \label{-1/2}
  p_1(a) \lesssim n^{-1/2}.
\end{equation}

This is sharp for $a_k = 1$: there are lots of cancelations in most
of the sign-sums $\sum \pm 1$. However, if $|a_j - a_k| \ge 1$ for
$k \ne j$, then the small ball probability is even smaller:
\begin{equation}                \label{-3/2}
  p_1(a) \lesssim n^{-3/2}.
\end{equation}
This was proved by Erd\"os and Moser \cite{E 65} for $p_0(a)$ and
with an extra $\log n$ factor, which was removed by S\'ark\"ozi and
Szemer\'edi \cite{SS}. H\'alasz \cite{H} proved this estimate for
$p_1(a)$ and generalized it to higher dimensions. Estimate
\eqref{-3/2} is sharp for $a_k = k$: there are still many
cancelations in most of the sign-sums $\sum \pm k$.

Tao and Vu \cite{TV} recently proposed a method to reduce the small
ball probability to an arbitrary polynomial order. They suggested to
look at the inverse problem and to study the following phenomenon:
\begin{quote}
{\em
  If the small ball probability $p_0(a)$ is large then
  the coefficient vector $a$ has a rich additive structure.
}
\end{quote}
Thus, the only reason for many cancelations in the sign-sums $\sum
\pm a_k$ is that most coefficients $a_k$ are arithmetically well
comparable. By removing this obstacle one can force the small ball
probability down to an arbitrary polynomial order:

\begin{theorem}[Tao, Vu \cite{TV}]          \label{t: Tao-Vu}
  Let $a_1,\ldots,a_n$ be integers, and let $A \ge 1$, $\e \in (0,1)$.
  Suppose for the random sign-sums one has
  $$
  p_0(a) \ge n^{-A}.
  $$
  Then all except $O_{A,\e}(n^{\e})$ coefficients $a_k$ are contained in
  the Minkowski sum of $O(A/\e)$
  arithmetic progressions of lengths $n^{O_{A,\e}(1)}$.
\end{theorem}

(Recall that the Minkowski sum of sets is defined as $U+V = \{ u+v :
\;  u \in U, \ v \in V \}$.)

\medskip

In  this paper we demonstrate that a similar, and even simpler,
phenomenon holds for real rather than integer numbers $a_k$, for the
small ball probabilities $p_\e(a)$ rather than the probability
$p_0(a)$ of exact values, and for general random sums \eqref{random
sum} rather than the random sign-sums.

We thus provide an essentially sharp solution to the
Littlewood-Offord problem for coefficients $a_k$ of equal order of
magnitude. We show that one can force the small ball probability
$p_{\e}(a)$ down to an {\em arbitrary function} of $n$, up to an
exponentially small order, which is best possible. We prove that:
\begin{quote}
{\em
  The coefficients of $a$ are essentially contained in one
  arithmetic progression of length $\lesssim p_{\e}(a)^{-1}$.
}
\end{quote}

By ``essentially'' we mean that for arbitrary $\a \in (0,1)$ and $\k
> c(\a)$ we can guarantee that all but $\k$ coefficients $a_k$ are
within $\a d$ from the elements of some arithmetic progression,
where $d$ is the gap between its elements. It is convenient to state
this result in terms of the essential least common denominator of
real numbers:

\begin{definition}[Essential LCD]
  Let $\a \in (0,1)$ and $\k \ge 0$.
  The {\em essential least common denominator}
  $D(a) = D_{\a,\k}(a)$ of a vector $a \in \R^n$ is defined as the
  infimum of $t > 0$ such that all except $\k$ coordinates of the vector
  $t a$ are of
  distance at most $\a$ from nonzero integers.
\end{definition}

For numbers $a_k = O(1)$, the essential LCD has an obvious
interpretation in terms of arithmetic progressions: all except $\k$
coefficients $a_k$ are within distance $\a/D(a) = O(\a)$ from the
elements of an arithmetic progression of length $O(D(a))$.

\begin{theorem}[Small Ball Probability]                \label{t: small ball intro}
  Let $\xi_1,\ldots,\xi_n$ be independent identically distributed
  centered random variables with variances at least $1$
  and third moments bounded by $B$.
  Let $a= (a_1, \ldots, a_n)$ be a vector of real coefficients such that,
  for some $K_1, K_2 > 0$ one has
  \begin{equation}                          \label{K1 K2}
     K_1 \le |a_k| \le K_2  \quad \text{for all $k$}.
  \end{equation}
  Let $\a \in (0,1)$ and $\k \in (0,n)$.
  Then for for every $\e \ge 0$ one has
  \[
    p_\e(a)
    \le \frac{C}{\sqrt{\k}} \Big( \e + \frac{1}{D_{\a,\k}(a)} \Big)
        + C e^{-c \a^2 \k},
  \]
  where $C, c > 0$ depend (polynomially) only on $B, K_1, K_2$.
\end{theorem}

A more precise version of this result is Theorem~\ref{t: small ball
precise} below.

\medskip

\remarks 1. By the definition, one always has $D_{\a,\k}(a) \gtrsim
1/K_2$ (e.g. with $\a =1/3$ and $\k = n/4$). Theorem~\ref{t: small
ball intro} thus yields $p_1(a) \lesssim n^{-1/2}$,
which agrees with Littlewood-Offord and Erd\"os inequality \eqref{-1/2}. 

2. Suppose the components of $a$ are uniformly spread between two
comparable values; say $a = (n, n+1, n+2, \ldots, 2n)$. 
Obviously, $D_{\a,\k}(a/n) \sim n$ (e.g. with $\a = 1/3$ and
$\k = n/4$). Theorem~\ref{t: small ball intro} thus yields
$p_1(a) = p_{1/n}(a/n) \lesssim n^{-3/2}$. This agrees
with Erd\"os-Moser inequality \eqref{-3/2}.

3. By making coefficients of $a$ more arithmetically incomparable,
such as by considering polynomial progressions, one can force the
small ball probability $p_\e(a)$ down to an arbitrarily small value, up
to an exponentially small order.

\medskip

One can restate Theorem \ref{t: small ball intro} as an inverse
Littlewood-Offord theorem:

\begin{corollary}[Inverse Littlewood-Offord Theorem]        \label{c: small ball}
  Let $a_1,\ldots,a_n$ be real numbers satisfying \eqref{K1 K2}
  and $\xi_1,\ldots,\xi_n$ be random variables as in
  Theorem~\ref{t: small ball intro}.
  Let $A \ge 1/2$, $\k \in (0,n)$ and $\e > 0$.
  Suppose for the random sums \eqref{random sum} one has
  $$
  p_\e(a) \ge n^{-A}.
  $$
  Then there exists an arithmetic progression of length $L = O(n^A \k^{-1/2})$
  and with gap between its elements $d \le 1$, and such that
  all except $\k$ coefficients $a_k$ are within distance
  $O(A \log(n) / \k)^{1/2} \cdot d$
  from the elements of the progression, provided that
  $\e \le 1/L$.
\end{corollary}

By Remark~1 above, the assumption $A \ge 1/2$ is optimal.

\smallskip

In contrast with Theorem~\ref{t: Tao-Vu}, Corollary~\ref{c: small
ball} guarantees an approximate, rather than exact, embedding of the
coefficients $a_1,\ldots,a_n$ into an arithmetic progression. On the
other hand, Corollary~\ref{c: small ball}: (a) applies for real
rather integer coefficients; (b) embeds into one arithmetic
progression rather than a  Minkowski sum of several progressions;
(c) provides a significantly sharper bound on the length of the
progression; (b) characterizes general small ball probabilities
$p_\e(a)$ rather than
  the probability of exact values $p_0(a)$;
(d) holds for general sums of i.i.d. random variables rather than
the
  random sign-sums.

\subsection{Outline of the argument}

We develop a general approach to the invertibility of random matrices.
Our main result, the Strong Invertibility Theorem~\ref{t: strong},
{\em reduces estimating the smallest singular value of random matrices
to estimating the largest singular value.}
Because the largest singular value is much more studied, this
immediately implies both our invertibility results stated above,
Theorems~\ref{t: 4} and ~\ref{t: subgaussian}.

The general approach to invertibility is developed in two stages. In
Section~\ref{s: median}  we present a ``soft'' and rather short
argument that leads to a weaker result. It yields
the Fourth Moment Theorem~\ref{t: 4} and also a weaker version of
the Subgaussian Theorem~\ref{t: subgaussian}
with $Cn^{-1/2}$ instead of the exponential term $c^n$
in \eqref{eq main}.

Our soft argument does not use any new estimates of the small ball probability.
To bound $\|Ax\|_2$ below for all vectors $x$ in the unit sphere, we
give two separate arguments for {\em compressible} vectors $x$,
whose norm is concentrated in a small number of coordinates, and for
{\em incompressible} vectors comprising the rest of the sphere.

For a compressible vector, the main contribution in the quantity
$\|Ax\|_2$ comes from the few (say, $n/10$) columns of $A$
corresponding to the biggest coordinates of $x$. This allows us to
replace $A$ by its $n \times n/10$ submatrix with the chosen
columns. Such rectangular random matrices are known to have big
smallest singular value (see e.g. \cite{LPRT}), which establishes a
nice lower bound on $\|Ax\|_2$ for all compressible vectors.

For the incompressible vectors, we show the invertibility differently.
Clearly, $s_n(A)$ is bounded above by the
distance from its $n$-th row vector $X_n$ to the span $H_n$ of the
others. We use a careful {\em average union} argument
(Lemma~\ref{l: via distance}) to show a reverse inequality
for $A$ restricted to the set of incompressible vectors.

Next, this distance can be bounded below as $\dist(X_n,H_n) \ge |\<
X^*, X_n \> |$, where $X^*$ is a unit normal of $H_n$. Since $X^*$
and $H_n$ are independent, the inner product $\< X^*, X_n \> $ can
be written as a sum of independent random variables of the form
\eqref{random sum}. This reduces the invertibility problem to the
Littlewood-Offord problem.

A useful small ball probability bound can be deduced from the
central limit theorem, by approximating the random sum \eqref{random
sum} with a Gaussian random variable for which the small ball
probability is easy to compute. With such bound, the argument above
yields a weaker version of the invertibility estimate \eqref{eq
main} with $Cn^{-1/2}$ instead of $c^n$.

This weaker estimate is a limitation of using the central limit theorem.
To prove the Strong Invertibility Theorem~\ref{t: strong},
and thus deduce the Subgaussian Theorem~\ref{t: subgaussian},
we will use the full strength of the Small Ball Probability
Theorem~\ref{t: small ball intro} instead.
This argument is presented in Section \ref{s: main proof}.

Our proof of  Theorem~\ref{t: small ball intro} starts with the
method developed by Hal\'asz \cite{H 75, H}. It allows us to bound
the small ball probability $p_\e(a)$ by a quantity of ergodic nature
-- the measure of the {\em recurrence set} of $a$. It indicates how
often a particle in $\R^n$ moving in the direction $a$ with unit
speed gets close to the points of the integer lattice. If this
happens often, then a density argument shows that the particle must
get close to two distinct lattice points over a short period of
time, say at times $t_1$ and $t_2$. It then follows that
$(t_2-t_1)a$ is close to an integer, which implies that the
essential LCD of $a$ is small. This argument is given in
Section~\ref{s: small ball}.

\subsection*{Acknowledgement.} The authors are grateful to the referee
for the careful reading of the manuscript and valuable suggestions.

\section{Preliminaries}
In the sequel $n$ denotes a sufficiently large integer, i.e. 
an integer bigger than a suitable absolute constant.
The standard inner product on $\R^n$ is denoted by $\< x, y \> $.
The $\ell_p$ norm on $\R^n$ is defined as
$\|x\|_p = (\sum_{k=1}^n |x_k|^p)^{1/p}$
for $0 < p < \infty$, and $\|x\|_\infty = \max_k |x_k|$.
The unit Euclidean ball and the sphere in $\R^n$ are denoted by
$B_2^n$ and $S^{n-1}$ respectively.
For a subset $\s \subseteq \{1,\ldots,n\}$, the orthogonal projection
onto $\R^\s$ in $\R^n$ is denoted by $P_\s$.

\medskip

The following observation will allow us to select
a nice subset of the coefficients $a_k$ when computing the small
ball probability.

\begin{lemma}[Restriction]                      \label{l: restriction}
  For any $a \in \R^n$, any $\s \subseteq \{1,\ldots,n\}$ and any $\e \ge 0$,
  we have
  $$
  p_\e(a) \le p_\e(P_\s a).
  $$
\end{lemma}

\begin{proof}
For fixed $v \in \R$ and for the random sum \eqref{random sum},
we write $S - v = S_\s - v_\s$, where
$S_\s := \sum_{k \in \s} a_k \xi_k$ and
$v_\s := v - \sum_{k \in \s^c} a_k \xi_k$.
We condition on a realization of $(\xi_k)_{k \in \s^c}$, and denote by $\P_\s$
the probability with respect to $(\xi_k)_{k \in \s}$.
Then a realization of $v_\s$ is fixed, so
$$
\P_\s (|S-v| \le \e) = \P_\s (|S_\s - v_\s| \le \e) \le p_\e(P_\s a).
$$
Taking the expectation of both sides with respect to $(\xi_k)_{k \in \s^c}$
completes the proof.
\end{proof}

The following tensorization lemma transfers one-dimensional
small ball probability estimates to the multidimensional case.
It is a minor variant of Lemma 4.4 of \cite{R}.

\begin{lemma}[Tensorization]   \label{l: tensorization}
  Let $\z_1, \ldots, \z_n$ be independent non-negative random
  variables, and let $K, \e_0 \ge 0$.

  (1) Assume that for each $k$
  $$
  \P (\z_k < \e) \le K \e
  \qquad \text{for all $\e \ge \e_0$}.
  $$
  Then
  $$
  \P \Big( \sum_{k=1}^n \z_k^2 < \e^2 n \Big) \le (C K \e)^n
  \qquad \text{for all $\e \ge \e_0$},
  $$
  where $C$ is an absolute constant.

  (2) Assume that there exist $\l>0$ and $\mu \in (0,1)$ such that for each $k$
  $$
  \P (\z_k < \l) \le \mu.
  $$
  Then there exist $\l_1>0$ and $\mu_1 \in (0,1)$
  that depend on $\l$ and $\mu$ only and such that
  $$
  \P \Big( \sum_{k=1}^n \z_k^2 < \l_1 n \Big) \le \mu_1^n.
  $$
\end{lemma}

We give a proof of the first part for completeness. The second
part is similar, cf. \cite{LPRT} proof of Proposition~3.4.

\begin{proof}
Let $\e \ge \e_0$. By Chebychev's inequality,
\begin{align}                              \label{tensorization product}
\P \Big( \sum_{k=1}^n \z_k^2 < \e^2 n \Big)
  = \P \Big( n - \frac{1}{\e^2} \sum_{k=1}^n \z_k^2 > 0 \Big)
  &\le \E \exp \Big( n - \frac{1}{\e^2} \sum_{k=1}^n \z_k^2 \Big) \notag \\
  &= e^n \prod_{k=1}^n \E \exp(-\z_k^2/\e^2).
\end{align}
By the distribution integral formula,
$$
\E \exp(-\z_k^2/\e^2) = \int_0^1 \P \big( \exp(-\z_k^2/\e^2) > s
\big) \; ds = \int_0^\infty 2 u e^{-u^2} \, \P(\z_k < \e u) \; du.
$$
For $u \in (0,1)$, we have $\P(\z_k < \e u) \le \P(\z_k < \e) \le
K \e$. This and the assumption of the lemma yields
$$
\E \exp(-\z_k^2/\e^2) \le \int_0^1 2 u e^{-u^2} K \e \; du
  + \int_1^\infty 2 u e^{-u^2} K \e u \; du
\le  C K \e.
$$

Putting this into \eqref{tensorization product} yields
$$
\P \Big( \sum_{k=1}^n \z_k^2 < \e^2 n \Big) \le e^n ( C K \e)^n.
$$
This completes the proof.
\end{proof}

\subsection{Largest singular value}

We recall some known bounds on the largest singular value of
random matrices under the fourth moment assumption and the subgaussian
moment assumption.
The following result is a partial case of a recent result of Latala.

\begin{theorem}[Largest singular value: fourth moment \cite{Lat}] \label{t: norm 4}
  Let $A$ be an $n \times n$ matrix whose entries are independent
  centered random variables with variances at least $1$ and
  fourth moments bounded by $B$. Then
  $$
  \E \|A\| \le C_1 n^{1/2}
  $$
  where $C_1 = C B^{1/4}$, and where $C$ is an absolute constant.
\end{theorem}

Under the stronger subgaussian moment assumption, a standard
observation shows that $\|A\| \sim n^{1/2}$ with exponentially
large probability (see e.g. \cite{DS} and \cite{LPRT}, Fact 2.4):

\begin{lemma}[Largest singular value: subgaussian]       \label{l: norm subgaussian}
  Let $A$ be an $n \times n$ matrix whose entries are independent
  centered random variables with variances at least $1$ and
  subgaussian moments bounded by $B$. Then
  $$
  \P ( \|A\| > C_1 n^{1/2} ) \le 2e^{-n},
  $$
  where $C_1$ depends only on $B$.
\end{lemma}

\subsection{Smallest singular value of rectangular matrices}

Estimates on the smallest singular value are known for
{\em rectangular} random matrices \cite{LPRT}.

\begin{proposition}[Smallest singular value of rectangular matrices]
            \label{p: rectangular}
  Let $G$ be an $n \times k$ matrix whose entries are independent
  centered random variables with variances at least $1$ and
  fourth moments bounded by $B$.
  Let $K \ge 1$.
  Then there exist $c_1, c_2 > 0$ and $\d_0 \in (0,1)$ that depend only
  on $B$ and $K$ such that if $k < \d_0 n$ then
  \begin{equation}                  \label{eq: rectangular}
    \P \big( \inf_{x \in S^{k-1}} \|Gx\|_2 \le c_1 n^{1/2}
        \text{ and } \norm{G} \le K n^{1/2} \big)
    \le e^{-c_2 n}.
  \end{equation}
\end{proposition}

Under the stronger subgaussian assumption, the condition
$\norm{G} \le K n^{1/2}$ can clearly be removed from \eqref{eq: rectangular}
by Lemma~\ref{l: norm subgaussian}.
This is not so under the fourth moment assumption. So here and later
in the paper, this condition will often appear in order to deduce
the Fourth Moment Theorem~\ref{t: 4}. The reader interested
only in the Subgaussian Theorem~\ref{t: subgaussian} can disregard this
condition.

A result stronger than Proposition~\ref{p: rectangular},
for the aspect ratio $\d_0$ arbitrarily close to 1,
follows by modifying the argument of \cite{LPRT}.
For completeness, we shall prove Proposition~\ref{p: rectangular}.
We start with the most general (but weakest possible) estimate
on the small ball probability.

\begin{lemma}   \label{l: LPRT}
  Let $\xi_1 \etc \xi_n$ be independent
  centered random variables with variances at least $1$ and
  fourth moments bounded by $B$.
  Then there exists $\mu \in (0,1)$ depending only on $B$,
  such that for every coefficient vector $a= (a_1,\ldots,a_n) \in S^{n-1}$
  the random sum $S = \sum_{k=1}^n a_k \xi_k$ satisfies
  $$
  \P(|S| < 1/2) \le \mu.
  $$
\end{lemma}

\begin{proof}
Let $\e_1 \etc \e_n$ be independent symmetric $\pm 1$ random variables,
which are independent of $\xi_1 \etc \xi_n$.
By the standard symmetrization inequality (see \cite{LT} Lemma~6.3),
$$
\E S^4 \le 16 \E \Big( \sum_{k=1}^n \e_k \xi_k a_k \Big)^4.
$$
We first condition on $\xi_1 \etc \xi_n$ and take the expectation with
respect to $\e_1 \etc \e_n$. Khinchine's inequality
(see e.g. \cite{LT} Lemma 4.1) and our assumptions on $\xi_k$ then yield
\begin{align*}
\E S^4
  &\le C \E \Big( \sum_{k=1}^n \xi_k^2 a_k^2 \Big)^2
  = C \E \sum_{k,j=1}^n \xi_k^2 \xi_j^2 a_k^2 a_j^2\\
  &\le C  \sum_{k,j=1}^n (\E \xi_k^4 )^{1/2} (\E \xi_j^4 )^{1/2} a_k^2 a_j^2
  \le CB \Big( \sum_{k=1}^n a_k^2 \Big)^2 = CB.
\end{align*}
The Paley--Zygmund inequality (see e.g. \cite{LPRT}, Lemma 3.5)
implies that for any $\l>0$
$$
\P (|S| > \l)
\ge \frac{(\E S^2 - \l^2)^2}{\E S^4}
     \ge \frac{(1-\l^2)^2}{CB}.
$$
To finish the proof, set $\l=1/2$.
\end{proof}

Combining Lemma \ref{l: LPRT} with the tensorization Lemma \ref{l:
tensorization}, we obtain the following invertibility
estimate for a fixed vector.

\begin{corollary}  \label{c: individual}
 Let $G$ be a matrix as in Proposition~\ref{p: rectangular}. Then
 there exist constants $\eta,\nu \in (0,1)$ depending only on $B$,
 such that for every $x \in S^{k-1}$
 \[
   \P ( \norm{Gx}_2  < \eta n^{1/2}) \le \nu^n.
 \]
\end{corollary}

\medskip

\begin{proof}[Proof of Proposition~\ref{p: rectangular}]
Let $\e>0$ to be chosen later. There exists an $\e$-net $\NN$
in $S^{k-1}$ (in the Euclidean norm) of cardinality
$|\NN| \le (3/\e)^k$ (see e.g. \cite{MS}).
Let $\eta$ and $\nu$ be the numbers in Corollary \ref{c: individual}.
Then by the union bound,
\begin{equation}  \label{on net}
  \P \left ( \exists x \in \NN : \ \norm{Gx}_2 < \eta n^{1/2} \right )
  \le (3/\e)^k \cdot \nu^n.
\end{equation}
Let $V$ be the event that $\norm{G} \le K n^{1/2}$ and
$\norm{Gy}_2 \le \frac{1}{2} \eta n^{1/2}$ for some point $y \in S^{k-1}$.
Assume that $V$ occurs, and choose a point $x \in \NN$ such that
$\norm{y-x}_2<\e$. Then
$$
\norm{Gx}_2 \le \norm{Gy}_2 + \norm{G} \cdot \norm{x-y}_2
 \le \frac{1}{2} \eta n^{1/2} + K n^{1/2} \cdot \e
 = \eta n^{1/2},
$$
if we set $\e=\eta/2K$. Hence, by \eqref{on net},
$$
    \P (V) \le \big( \nu \cdot  \left ( 3/\e \right
    )^{k/n} \big)^n
    \le e^{-c_2 n},
$$
if we assume that $k/n \le \d_0$ for an appropriately chosen $\d_0 < 1$.
This completes the proof.
\end{proof}

\subsection{The small ball probability via the central limit theorem}

The central limit theorem can be used to estimate
the small ball probability, as observed in \cite{LPRT}.
Specifically, one can use the Berry-Ess\'een version of the central limit
theorem (see \cite{Str}, Section~2.1):

\begin{theorem}[Berry-Ess\'een CLT]                             \label{CLT}
  Let $\z_1,\ldots,\z_n$ be independent centered random variables
  with finite third moments, and let
  $\s^2 := \sum_{k=1}^n \E |\z_k|^2$.
  Consider a standard normal random variable $g$.
  Then for every $t > 0$:
  \begin{equation}                          \label{eq CLT}
    \Big| \P \Big( \frac{1}{\s} \sum_{k=1}^n \z_k \le  t \Big) - \P(g \le t) \Big|
    \le C \s^{-3} \sum_{k=1}^n \E |\z_k|^3,
  \end{equation}
  where $C$ is an absolute constant.
\end{theorem}

The following corollary is essentially given in \cite{LPRT}.
We shall include a proof for the reader's convenience.  

\begin{corollary}[Small ball probability via CLT]         \label{SBP via CLT}
  Let $\xi_1 \etc \xi_n$ be independent centered random variables
  with variances at least $1$ and third moments bounded by $B$.
  Then for every $a \in \R^n$ and every $\e \ge 0$, one has
  \[
    p_\e(a)
    \le \sqrt{\frac{2}{\pi}} \frac{\e}{\|a\|_2}
      + C_1 B \Big( \frac{\|a\|_3}{\|a\|_2} \Big)^3,
  \]
  where $C_1$ is an absolute constant.
\end{corollary}

\begin{proof}
We shall use
Theorem~\ref{CLT} for $\z_k = a_k \xi_k$. There, $\s \ge \norm{a}_2$
and $\sum_{k=1}^n \E |\z_k|^3 \le B \norm{a}_3^3$. Thus for every $u
\in \R$ we have
\begin{equation}                                \label{SBP ut}
  \P \Big( \Big| \frac{1}{\|a\|_2} \sum_{k=1}^n a_k \xi_k - u \Big| \le t \Big)
  \le \P ( |g-u| \le t ) + 2 C B \Big( \frac{\|a\|_3}{\|a\|_2} \Big)^3.
\end{equation}
Since the density of the standard normal random variable $g$ is uniformly
bounded by $1/\sqrt{2\pi}$, we have
$$
\P ( |g-u| \le t ) \le \frac{2t}{\sqrt{2\pi}} = \sqrt{\frac{2}{\pi}} \; t.
$$
With $u = \frac{v}{\|a\|_2}$ and $t = \frac{\e}{\|a\|_2}$,
the left hand side of \eqref{SBP ut} equals $\P (|S - v| \le \e)$,
which completes the proof with $C_1 = 2C$.
\end{proof}

As an immediate corollary, we get:

\begin{corollary}[Small ball probability for big $\e$]       \label{SBP big e}
  Let $\xi_1 \etc \xi_n$ be independent centered random variables
  with variances at least $1$ and third moments bounded by $B$.
  Assume that a coefficient vector $a$ satisfies \eqref{K1 K2}.
  Then for every $\e \ge 0$ one has
  $$
  p_\e(a) \le \frac{C_2}{\sqrt{n}} \; \big( \e/K_1 + B (K_2/K_1)^3 \big),
  $$
  where $C_2$ is an absolute constant.
\end{corollary}

\section{Invertibility of random matrices: soft approach}
\label{s: median}

In this section, we develop a soft approach to the invertibility
of random matrices. Instead of using the new estimates on the small
ball probability, we will rely on the central limit theorem
(Corollary~\ref{SBP big e}). This approach will yield
a weaker bound, with polynomial rather than exponential term
for the singularity probability. In Section~\ref{s: main proof} we
shall improve upon the weak point of this argument, so
the Small Ball Probability Theorem~\ref{t: small ball intro}
will be used instead.

\begin{theorem}[Weak invertibility]                 \label{t: weak}
  Let $A$ be an $n \times n$ matrix whose entries are independent
  random variables with variances at least $1$ and fourth moments
  bounded by $B$.
  Let $K \ge 1$. Then for every $\e \ge 0$ one has
  \begin{equation}              \label{eq: weak}
  \P \big( s_n(A) \le \e n^{-1/2} \big)
    \le C\e + Cn^{-1/2} + \P (\|A\| > K n^{1/2}),
  \end{equation}
  where $C$ depends (polynomially) only on $B$ and $K$.
\end{theorem}

To make this bound useful, we recall that the last term in \eqref{eq: weak}
can be bounded using Theorem~\ref{t: norm 4} under the fourth
moment assumption and by Lemma~\ref{l: norm subgaussian} under the subgaussian
assumption. In particular, this proves Fourth Moment Theorem~\ref{t: 4}:

\begin{proof}[Proof of the Fourth Moment Theorem~\ref{t: 4}]
Let $\d > 0$. By Theorem~\ref{t: norm 4} and using Chebychev's inequality,
we have
$$
\P \Big( \|A\| > \frac{3C_1}{\d} n^{1/2} \Big) < \d/3.
$$
Then setting $K = 3C_1/\d$, $\e = \d/3C$ and $n_0 = (3C/\d)^2$,
we make each of the three terms in the right hand side of \eqref{eq: weak}
bounded by $\d/3$. This completes the proof.
\end{proof}

\begin{remark}
  Theorem~\ref{t: weak} in combination with Lemma~\ref{l: norm subgaussian}
  yields a weaker version of the Subgaussian
  Theorem~\ref{t: subgaussian}, with $Cn^{-1/2}$ instead of $c^n$.
\end{remark}

\subsection{Decomposition of the sphere}            \label{s: decomposition}

To prove Theorem~\ref{t: weak}, we
shall partition the unit sphere $S^{n-1}$ into the two sets of
{\em compressible} and {\em incompressible} vectors, and will show
the invertibility of $A$ on each set separately.

\begin{definition}[Compressible and incompressible vectors]   \label{d: compressible}
  Let $\d, \rho \in (0,1)$.
  A vector $x \in \R^n$ is called {\em sparse} if
  $|\supp(x)| \le \d n$.
  A vector $x \in S^{n-1}$ is called {\em compressible} if $x$
  is within Euclidean distance $\rho$ from the set of
  all sparse vectors.
  A vector $x \in S^{n-1}$ is called {\em incompressible}
  if it is not compressible.
  The sets of sparse, compressible and incompressible vectors
  will be denoted by $\Sparse = \Sparse(\d)$, $\Comp = \Comp(\d,\rho)$
  and $\Incomp = \Incomp(\d,\rho)$ respectively.
\end{definition}

\begin{remarks}
  1. Here we borrow the terminology from the
  signal processing and the sparse approximation theory.
  Efficient compression of many real-life signals, such as
  images and sound, relies on the assumption that their
  coefficients (Fourier, wavelet, frame etc.) decay in a fast way.
  Essential information
  about the signal is thus contained in few most significant
  coefficients, which can be stored in small space
  (see \cite{Do, CT}).
  Such coefficient vector is close to a sparse
  vector, and is thus compressible in the sense of our definition.

  2. Sets similar to those of compressible and incompressible
  vectors were previously used for the invertibility problem
  in \cite{LPRT} and \cite{R}.

  3. In our argument, the parameters $\d, \rho$ will be chosen as small
  constants that depend only on $B$ and $K$.
\end{remarks}

\medskip

Using the decomposition of the sphere $S^{n-1} = \Comp \cup \Incomp$,
we break the invertibility problem into two subproblems,
for compressible and incompressible vectors:
\begin{multline}                        \label{two terms}
\P \big( s_n(A) \le \e n^{-1/2} \text{ and } \|A\| \le K n^{1/2} \big) \\
  \le \P \big( \inf_{x \in \Comp(\d,\rho)} \|A x\|_2 \le \e n^{-1/2}
        \text{ and } \|A\| \le K n^{1/2} \big) \\
    + \P \big( \inf_{x \in \Incomp(\d,\rho)} \|A x\|_2 \le \e n^{-1/2}
        \text{ and } \|A\| \le K n^{1/2} \big).
\end{multline}

The compressible vectors are close to a coordinate subspace of a
small dimension $\d n$. The restriction of our random matrix
$A$ onto such a subspace is a random {\em rectangular} $n \times \d n$
matrix. Such matrices are well invertible with exponentially high probability
(see Proposition~\ref{p: rectangular}).
By taking the union bound over all coordinate subspaces,
we will deduce the invertibility of the random matrix on the
set of compressible vectors.

Showing the invertibility on the set of incompressible vectors is generally
harder, for this set is bigger in some sense.
By a careful {\em average union} argument, we shall reduce
the problem to a small ball probability estimate.

\subsection{Invertibility for the compressible vectors}

On the set of compressible vectors, a much stronger invertibility holds
than we need in \eqref{two terms}:

\begin{lemma}[Invertibility for compressible vectors]       \label{l: compressible}
  Let $A$ be a random matrix as in Theorem~\ref{t: weak},
  and let $K \ge 1$.
  Then there exist $\d, \rho, c_3, c_4 > 0$ that depend only
  on $B$ and $K$, and such that
  $$
  \P \big( \inf_{x \in \Comp(\d,\rho)} \|A x\|_2 \le c_3 n^{1/2}
    \text{ and } \|A\| \le K n^{1/2} \big)
  \le e^{-c_4 n}.
  $$
\end{lemma}

\begin{remark}
  The bound in Lemma~\ref{l: compressible} is much stronger than
  we need in \eqref{two terms}. Indeed, by choosing the constant $C$
  in Theorem~\ref{t: weak} large enough, we can assume that
  $n > 1/c_3$ and $\e < 1$. Then the value $c_3 n^{1/2}$
  in Lemma~\ref{l: compressible} is bigger than $\e n^{-1/2}$
  in \eqref{two terms}.
\end{remark}

\begin{proof}
We first prove a similar invertibility estimate for the sparse
vectors. To this end, we can assume that $\d_0 < 1/2$ in
Proposition~\ref{p: rectangular}. We use this result with $k = \d n$
and take the union bound over all
$\lceil \d n \rceil$-element subsets $\s$ of $\{1,\ldots,n\}$:
\begin{align}                                                   \label{sparse}
\P
  &\big( \inf_{x \in \Sparse(\d), \, \|x\|_2 = 1} \|Ax\|_2 \le c_1 n^{1/2}
    \text{ and } \|A\| \le K n^{1/2} \big) \\
  &= \P \big( \exists \s, \; |\s| = \lceil \d n \rceil :
    \inf_{x \in \R^\s,\; \|x\|_2 = 1} \|Ax\|_2 \le c_1 n^{1/2}
    \text{ and } \|A\| \le K n^{1/2} \big) \notag\\
  &\le \binom{n}{\lceil \d n \rceil} \, e^{-c_2 n}
  \le \exp(4 e \d \log(e/\d) n - c_2 n)
  \le e^{-c_2 n / 2} \notag
\end{align}
with an appropriate choice of $\d < \d_0$,
which depends only on $c_2$ (which in turn depends only on $B$ and $K$).

Now we deduce the invertibility estimate for the compressible vectors.
Let $c_3 > 0$ and $\rho \in (0,1/2)$ to be chosen later.
We need to bound the event $V$ that $\|Ax\|_2 \le c_3 n^{1/2}$
for some vector $x \in \Comp(\d,\rho)$ and $\|A\| \le K n^{1/2}$.
Assume $V$ occurs.
Every such vector $x$ can be written as a sum
$x = y + z$, where $y \in \Sparse(\d)$ and $\|z\|_2 \le \rho$.
Thus $\|y\|_2 \ge 1-\rho \ge 1/2$, and
$$
\|Ay\|_2 \le \|Ax\|_2 + \|A\|\|z\|_2
\le c_3 n^{1/2} + \rho K n^{1/2}.
$$
We choose $c_3 := c_1/4$ and $\rho := c_1/4K$ so that
$\|Ay\|_2 \le \frac{1}{2} c_1 n^{1/2}$. Since $\|y\|_2 \ge 1/2$, we
have found a unit vector $u \in \Sparse(\d)$ such that
$\|Au\|_2 \le c_1 n^{1/2}$ (choose $u = y/\|y\|_2$).
This shows that the event $V$ implies the event in \eqref{sparse}, so we have
$\P(V) \le e^{-c_2 n / 2}$.
This completes the proof.
\end{proof}

\subsection{Invertibility for the incompressible vectors via distance}

For the incompressible vectors, we shall reduce the invertibility problem
to a lower bound on the distance between a random vector and a random hyperplane.

We first show that incompressible vectors are well spread in the
sense that they have many coordinates of the order $n^{-1/2}$.

\begin{lemma}[Incompressible vectors are spread]    \label{l: spread}
  Let $x \in \Incomp(\d,\rho)$.
  Then there exists a set $\s \subseteq \{1, \ldots, n\}$
  of cardinality $|\s| \ge \frac{1}{2} \rho^2 \d n$ and such that
  $$
  \frac{\rho}{\sqrt{2n}} \le |x_k| \le \frac{1}{\sqrt{\d n}}
  \qquad \text{for all $k \in \s$.}
  $$
\end{lemma}

\begin{proof}
Consider the subsets of $\{1,\ldots,n\}$ defined as
$$
\s_1 := \{ k : \; |x_k| \le \frac{1}{\sqrt{\d n}} \}, \qquad
\s_2 := \{ k : \; |x_k| \ge \frac{\rho}{\sqrt{2n}} \},
$$
and put $\s := \s_1 \cap \s_2$.

By Chebychev's inequality, $|\s_1^c| \le \d n$.
Then $y := P_{\s_1^c} x \in \Sparse(\d n)$, so the
incompressibility of $x$ implies that
$\|P_{\s_1} x\|_2 = \|x - y\|_2 > \rho$.
By the definition of $\s_2$, we have
$\|P_{\s_2^c} x\|_2^2 \le n \cdot \frac{\rho^2}{2n} = \rho^2/2$.
Hence
\begin{equation}                            \label{proj small}
\|P_\s x\|_2^2 \ge \|P_{\s_1} x\|_2^2 - \|P_{\s_2^c} x\|_2^2
\ge \rho^2/2.
\end{equation}
On the other hand, by the definition of $\s_1 \supseteq \s$,
\begin{equation}                            \label{proj big}
\|P_\s x\|_2^2 \le \|P_\s x\|_\infty^2 \cdot |\s|
\le \frac{1}{\d n} \cdot |\s|.
\end{equation}
It follows from \eqref{proj small} and \eqref{proj big} that
$|\s| \ge \frac{1}{2} \rho^2 \d n$.
\end{proof}

\begin{lemma}[Invertibility via distance]  \label{l: via distance}
  Let $A$ be any random matrix.
  Let $X_1,\ldots,X_n$ denote the column vectors of $A$, and let $H_k$ denote
  the span of all column vectors except the $k$-th.
  Then for every $\d,\rho \in (0,1)$ and every $\e > 0$, one has
  \begin{equation}                              \label{eq: via distance}
    \P \big( \inf_{x \in \Incomp(\d,\rho)} \|A x\|_2 < \e \rho n^{-1/2} \big)
    \le \frac{1}{\d n} \sum_{k=1}^n
          \P \big( \dist( X_k, H_k) < \e \big).
  \end{equation}
\end{lemma}

\begin{remark}
  The main point of this bound is the {\em average}, rather than the maximum,
  of the distances in the right hand side of \eqref{eq: via distance}.
  This will allow us to avoid estimating the union of $n$ events
  and thus bypass a loss of the $n$ factor in the invertibility theorem.
\end{remark}

\begin{proof}
Let $x \in \Incomp(\d,\rho)$.
Writing $A x = \sum_{k=1}^n x_k X_k$, we have
\begin{align}    \label{norm Q_j}
\|Ax\|_2
  &\ge \max_{k=1 \etc n} \dist (Ax, H_k) \notag \\
  &= \max_{k=1 \etc n} \dist (x_k X_k, H_k)
  = \max_{k=1 \etc n} |x_k| \, \dist (X_k, H_k).
\end{align}
Denote
$$
p_k := \P \big( \dist( X_k, H_k) < \e \big).
$$
Then
$$
\E \big| \{ k :\, \dist(X_k, H_k) < \e \} \big|
  = \sum_{k=1}^n p_k.
$$
Denote by $U$ the event that the set
$\s_1 := \{ k :\, \dist(X_k, H_k) \ge \e \}$
contains more than $(1-\d)n$ elements.
Then by Chebychev's inequality,
$$
\P(U^c) \le \frac{1}{\d n} \sum_{k=1}^n p_k.
$$
On the other hand, for every incompressible vector $x$, the set
$\s_2(x) := \{ k :\; |x_k| \ge \rho n^{-1/2} \}$
contains at least $\d n$ elements.
(Otherwise, since $\|P_{\s_2(x)^c} x\|_2 \le \rho$,
we would have $\|x - y\|_2 \le \rho$ for the sparse vector
$y := P_{\s_2(x)} x$, which would contradict the
incompressibility of $x$).

Assume that the event $U$ occurs.
Fix any incompressible vector $x$.
Then $|\s_1| + |\s_2(x)| > (1-\d)n + \d n > n$, so
the sets $\s_1$ and $\s_2(x)$ have nonempty intersection.
Let $k \in \s_1 \cap \s_2(x)$. Then by \eqref{norm Q_j}
and by the definitions of the sets $\s_1$ and $\s_2(x)$, we have
$$
\|Ax\|_2 \ge  |x_k| \, \dist (X_k, H_k)
  \ge \rho n^{-1/2} \cdot \e.
$$
Summarizing, we have shown that
$$
\P \big( \inf_{x \in \Incomp(\d,\rho)} \|A x\|_2 < \e \rho n^{-1/2} \big)
  \le \P(U^c)
  \le \frac{1}{\d n} \sum_{k=1}^n p_k.
$$
This completes the proof.
\end{proof}

\subsection{Distance via the small ball probability}        \label{s: dist via SBP}

Lemma~\ref{l: via distance} reduces the invertibility problem
to a lower bound on the distance between a random vector and a
random hyperplane. Now we reduce bounding the distance to a
small ball probability estimate.

Let $X_1,\ldots,X_n$ be the column vectors of $A$. These are
independent random vectors in $\R^n$. Consider the subspace $H_n =
\Span(X_1,\ldots,X_{n-1})$. Our goal is to bound the distance
between the random vector $X_n$ and the random subspace $H_n$.

To this end, let $X^*$ be any unit vector orthogonal to $X_1,\ldots, X_{n-1}$.
We call it a {\em random normal}. We can choose $X^*$ so that it is a
random vector that depends only on $X_1,\ldots,X_{n-1}$
and is independent of $X_n$.

We clearly have
\begin{equation}                                    \label{dist}
  \dist(X_n,H_n) \ge |\< X^*, X_n\> |.
\end{equation}
Since the vectors $X^* =: (a_1,\ldots,a_n)$ and
$X_n =: (\xi_1,\ldots,\xi_n)$ are independent, we should be able to use
the small ball probability estimates, such as Corollary~\ref{SBP big e},
to deduce a lower bound on the magnitude of
$$
\< X^*, X_n\> = \sum_{k=1}^n a_k \xi_k.
$$
To this end, we first need to check that the coefficients
of the vector $X^*$ are well spread.

\begin{lemma}[Random normal is incompressible]      \label{normal inc}
  Let $\d,\rho,c_4 > 0$ be as in Lemma~\ref{l: compressible}.
  Then
  $$
  \P \big( X^* \in Comp(\d,\rho) \text{ and } \|A\| \le K n^{1/2} \big)
    \le e^{-c_4 n}.
  $$
\end{lemma}

\begin{proof}
Let $A'$ be the $(n-1) \times n$ random matrix with rows
$X_1,\ldots,X_{n-1}$, i.e. the submatrix of $A^T$ obtained by removing
the last row. By the definition of the random normal,
\begin{equation}                                \label{kernel}
  A' X^* = 0.
\end{equation}
Therefore, if $X^* \in Comp(\d,\rho)$ then
$\inf_{x \in \Comp(\d,\rho)} \|A'x\|_2 = 0$.
By replacing $n$ with $n-1$, one can easilty check that the proof 
Lemma~\ref{l: compressible} remains valid for $A'$ as well as for $A$;
note also that $\|A'\| \le \|A\|$.
This completes the proof.
\end{proof}

Now we recall our small ball probability estimate, Corollary~\ref{SBP big e},
in a form useful for the incompressible vectors:

\begin{lemma}[Small ball probability for incompressible vectors]
                            \label{SBP incompressible}
  Let $\xi_1,\ldots,\xi_n$ be random variables as in Corollary~\ref{SBP big e}.
  Let $\d,\rho \in (0,1)$, and
  consider a coefficient vector $a \in \Incomp(\d,\rho)$.
  Then for every $\e \ge 0$ one has
  $$
  p_\e(a) \le C_5 (\e + B n^{-1/2}),
  $$
  where $C_5$ depends (polynomially) only on $\d$ and $\rho$.
\end{lemma}

\begin{proof}
Let $\s$ denote the set of the spread coefficients of $a$ constructed
in Lemma~\ref{l: spread}. Then $|\s| \ge \frac{1}{2} \rho^2 \d n$, and the vector
$b := n^{1/2} P_\s a$ satisfies $K_1 \le |b_k| \le K_2$ for all $k \in \s$,
where $K_1 = \rho/\sqrt{2}$ and $K_2 = 1/\sqrt{\d}$.
By Restriction Lemma~\ref{l: restriction} and
Corollary~\ref{SBP big e}, we have
$$
p_\e(a) = p_{n^{1/2} \e}(n^{1/2} a)
\le p_{n^{1/2} \e}(b)
\le C_5 \; (\e + B n^{-1/2}).
$$
This completes the proof.
\end{proof}

Lemmae~\ref{SBP incompressible} and \ref{normal inc} imply the
desired distance bound:

\begin{lemma}[Weak Distance Bound]        \label{l: dist}
  Let $A$ be a random matrix as in Theorem~\ref{t: weak}.
  Let $X_1,\ldots,X_n$ denote its column vectors, and consider
  the subspace $H_n = \Span(X_1,\ldots,X_{n-1})$.
  Let $K \ge 1$. Then for every $\e \ge 0$, one has
  $$
  \P \big( \dist(X_n, H_n) < \e \text{ and } \|A\| \le K n^{1/2} \big)
  \le C_6(\e + n^{-1/2}),
  $$
  where $C_6$ depends only on $B$ and $K$.
\end{lemma}

\begin{remark}
  In Theorem~\ref{t: distance} below, we shall improve
  this distance bound by reducing the polynomial term $n^{-1/2}$
  by the exponential term $e^{-cn}$.
\end{remark}

\begin{proof}
We condition upon a realization of the random vectors
$X_1,\ldots,X_{n-1}$. This fixes realizations of the subspace $H_n$
and the random normal $X^*$. Recall that $X_n$ is independent of
$X^*$. We denote the probability with respect to $X_n$ by $\P_n$,
and the expectation with respect to $X_1,\ldots,X_{n-1}$ by
$\E_{1,\ldots,n-1}$. Then
\begin{multline}                           \label{c}
\P \big( |\< X^*, X_n\> | < \e \text{ and } \|A\| \le K n^{1/2} \big) \\
  \le \E_{1,\ldots,n-1} \P_n \big( |\< X^*, X_n\> | < \e
      \text{ and } X^* \in \Incomp(\d,\rho) \big) \\
    + \P \big( X^* \in \Comp(\d,\rho) \text{ and } \|A\| \le K n^{1/2} \big).
\end{multline}

Fix $\d,\rho > 0$ so that the conclusion of Lemma~\ref{normal inc}
holds. This bounds the last term in the right hand side of \eqref{c}
by $e^{-c_4 n}$.
Furthermore, by Lemma~\ref{SBP incompressible}, for any fixed realization
of $X_1,\ldots,X_n$ such that $X^* \in \Incomp(\d,\rho)$ we have
$$
\P_n \big( |\< X^*, X_n\> | < \e \big)
\le C_5' \; (\e + n^{-1/2}),
$$
where $C_5'$ depends only on $B$ and $K$.
It follows that
$$
\P \big( |\< X^*, X_n\> | < \e \text{ and } \|A\| \le K n^{1/2} \big)
\le C_5' \; (\e + n^{-1/2}) + e^{-c_4 n}.
$$
By \eqref{dist}, the proof is complete.
\end{proof}

Combining Lemma~\ref{l: via distance} and Lemma~\ref{l: dist},
we have shown the invertibility of a random matrix on the
set of incompressible vectors:

\begin{lemma}[Invertibility for incompressible vectors]       \label{l: incompressible}
  Let $A$ be a random matrix as in Theorem~\ref{t: weak}.
  Let $K \ge 1$ and $\d,\rho \in (0,1)$.
  Then for every $\e \ge 0$, one has
  $$
  \P \big( \inf_{x \in \Incomp(\d,\rho)} \|A x\|_2 \le \e \rho n^{-1/2} \big)
  \le \frac{C_7}{\d} (\e + n^{-1/2}) + \P (\|A\| > K n^{1/2}),
  $$
  where $C_7$ depends only on $B$ and $K$.
\end{lemma}

\subsection{Invertibility on the whole sphere}

The Weak Invertibility Theorem~\ref{t: weak}
now follows from the decomposition of the sphere \eqref{two terms}
into compressible and incompressible vectors,
and from the invertibility on each of the two parts
established in Lemma~\ref{l: compressible} (see the remark below it)
and Lemma~\ref{l: incompressible}
(used for $\d, \rho$ as in Lemma~\ref{l: compressible} and for $\e/\rho$
rather than $\e$). \qed

\section{Small ball probability}  \label{s: small ball}

In this section, we prove the following more precise version
of Theorem~\ref{t: small ball intro}.

\begin{theorem}[Small Ball Probability]               \label{t: small ball precise}
  Let $\xi$ be a centered random variable with variance at least $1$
  and with the third moment bounded by $B$.
  Consider independent copies $\xi_1,\ldots,\xi_n$ of $\xi$.
  Let $a= (a_1, \ldots, a_n)$ be a coefficient vector
  and let $K \ge 1$ be such that
  \begin{equation}                          \label{12}
    1 \le |a_k| \le K  \quad \text{for all $k$}.
  \end{equation}
  Let $0 < \a < 1/6K$ and $0 < \k < n$.
  Then for every $\e \ge 0$ one has
  \[
    p_\e(a)
    \le \frac{C B K^3}{\sqrt{\k}} \Big( \e + \frac{1}{D_{2\a,2\k}(a)} \Big)
        + C \exp \Big( -\frac{c \a^2 \k}{B^2} \Big),
  \]
  where $C, c > 0$ are absolute constants.
\end{theorem}

\begin{remark}
  1. This result clearly implies Theorem~\ref{t: small ball intro}.
  (Indeed, in Theorem~\ref{t: small ball intro} one can assume that
  $K_1 = 1$ by rescaling the coefficients $a_k$, and that $\a < 1/6K_2$
  by considering $\a/6K_2$ instead of $\a$.)

  2. Since the definition of $p_\e(a)$ includes shifts,
  Theorem~\ref{t: small ball precise} holds also for the shifted
  random variables $\xi_j'=\xi_j+t_j$ for any real numbers $t_1 \etc
  t_n$.
\end{remark}

The approach based on the central limit theorem establishes
Theorem~\ref{t: small ball precise} for the values
of $\e$ of constant order and above. Indeed, for $\e > \e_0 > 0$,
Corollary~\ref{SBP big e} yields
$$
p_\e(a) \le \frac{C_2' B K^3}{\sqrt{n}} \; \e
$$
where $C_2'$ depends only on $\e_0$.

For $\e$ below the constant order, this bound can not hold
without any additional information about the coefficient vector $a$.
Indeed, if all $a_k = 1$ then random sign-sums satisfy
$p_0(a) \ge \P (S = 0) \sim n^{-1/2}$.

We thus need to develop a tool sharper that the central limit theorem
to handle smaller $\e$.
Our new method uses the approach of Hal\'asz \cite{H 75, H},
which was also used in \cite{R}.

\subsection{Initial reductions, symmetrization, truncation}

Throughout the proof, absolute constants
will be denoted by $C$, $c$, $c_1$, \ldots
The particular value of each constant
can be different in different instances.

As explained above, we can assume in the sequel that
$\e$ is below a constant, such as
\begin{equation}                            \label{delta small}
  \e < \pi/4.
\end{equation}

We can also assume that $\k < n/2$ and that $a_k \ge 1$ by replacing,
if necessary, $\xi_k$ by $-\xi_k$.

We shall symmetrize the random variables $\xi_k$ and remove any
small values they can possibly take.
For many random variables, such as random $\pm 1$, this step is
not needed.

Let $\xi'$ be an independent copy of $\xi$ and define the random variable
$\z := |\xi - \xi'|$. Then
$$
\E \z^2 = 2 \E |\xi|^2 \ge 2
\quad \text{and} \quad
\E \z^3 \le 8 \E |\xi|^3 \le 8B.
$$
The Paley-Zygmund inequality (see e.g. \cite{LPRT}, Lemma 3.5)
implies that
\begin{equation}                        \label{pz}
  \P (\z > 1) \ge \frac{(\E \z^2 - 1)^3}{(\E \z^3)^2}
  \ge \frac{1}{64 B^2} =: \b.
\end{equation}
Denote by $\bar{\z}$ the random variable $\z$ conditioned on $\z > 1$.
Formally, $\bar{\z}$ is a random variable such that for every
measurable function $f$ one has
$$
\E f(\bar{\z}) = \frac{1}{\P (\z > 1)} \;
  \E f(\z) \one_{\{\z > 1\}}.
$$
It then follows by \eqref{pz} that for every measurable non-negative
function $f$, one has
\begin{equation}                            \label{zeta zetabar}
\E f(\z) \ge \b \; \E f(\bar{\z}).
\end{equation}

\subsection{Small ball probability via characteristic functions}

An inequality of Ess\'een (\cite{Ess}, see also \cite{H}), bounds
the small ball probability of a random variable $S$ by the $L_1$
norm of its characteristic function
$$
\phi(t) = \phi_S(t) = \E \exp(i S t).
$$

\begin{lemma}[Ess\'een's Inequality]                           \label{Esseen}
  For every random variable $S$ and for every $\e > 0$, one has
  $$
  \sup_{v \in \R} \P (|S - v| \le \e)
  \le C \int_{-\pi/2}^{\pi/2} |\phi(t/\e)| \; dt,
  $$
  where $C$ is an absolute constant.
\end{lemma}

We want to use Ess\'een's Inequality for the random sum
$S = \sum_{k=1}^n a_k \xi_k$.
The characteristic function of $a_k \xi_k$ is
$$
\phi_k(t) := \E \exp(i a_k \xi_k t) = \E \exp(i a_k \xi t),
$$
so the characteristic function of $S$ is then
$$
\phi(t) = \prod_{k=1}^n \phi_k(t).
$$
To estimate the integral in Ess\'een's Lemma~\ref{Esseen},
we first observe that
$$
|\phi_k(t)|^2 = \E \cos(a_k \z t).
$$
Using the inequality $|x| \le \exp (-\frac{1}{2} (1-x^2))$ valid for
all $x$, we then obtain
\begin{align*}
|\phi(t)|
  &\le \prod_{k=1}^n \exp \Big( - \frac{1}{2} (1 - |\phi_k(t)|^2) \Big)\\
  &= \exp \Big( -\E \sum_{k=1}^n \frac{1}{2} (1 - \cos(a_k \z t)) \Big)
  = \exp \big( -\E f(\z t) \big),
\end{align*}
where
$$
f(t) := \sum_{k=1}^n \sin^2 \big( \frac{1}{2} a_k  t \big).
$$
Hence by \eqref{zeta zetabar}, we have
$$
|\phi(t)| \le \exp \big( - \b \; \E f(\bar{\z} t) \big).
$$
Then by Ess\'een's Lemma~\ref{Esseen} and using Jensen's inequality,
we estimate the small ball probability as
\begin{align}                           \label{SBP via int}
p_\e(a)
  \le C \int_{-\pi/2}^{\pi/2} |\phi(t/\e)| \; dt
  &\le C \int_{-\pi/2}^{\pi/2} \exp \big( - \b \; \E f(\bar{\z} t / \e) \big) \; dt \nonumber\\
  &\le C \E \; \int_{-\pi/2}^{\pi/2} \exp \big( - \b f(\bar{\z} t / \e) \big) \; dt \nonumber\\
  &\le C \sup_{z \ge 1} \;
    \int_{-\pi/2}^{\pi/2} \exp \big( - \b f(zt/\e) \big) \; dt.
\end{align}

Fix $z \ge 1$. First we estimate the maximum
$$
M := \max_{|t| \le \pi/2} f(zt/\e)
= \max_{|t| \le \pi/2} \sum_{k=1}^n \sin^2 (a_k z t / 2 \e).
$$

\begin{lemma}                               \label{M}
  We have
  $$
  \frac{n}{4}  \le M \le n.
  $$
\end{lemma}

\begin{proof}
The upper bound is trivial. For the lower bound, we estimate the
maximum by the average:
\[
  M \ge \frac{1}{\pi}
     \int_{-\pi/2}^{\pi/2} f(z t/\e) \,
     dt
   = \frac{1}{2} \sum_{k=1}^n
       \left ( 1 - \frac{\sin(\pi a_k z/2 \e)}{\pi a_k z/2 \e} \right ).
\]
By our assumptions, $a_k \ge 1$, $z \ge 1$ and $\e < \pi/4$. Hence
$\pi a_k z/2 \e \ge 2$, so
\[
  M \ge \frac{n}{2} \; \inf_{t \ge  2 } \Big( 1 - \frac{\sin t}{t} \Big)
  \ge \frac{n}{4}.
\]
This completes the proof.
\end{proof}

\medskip

Now we consider the level sets of $f$, defined for $m, r \ge 0$ as
$$
T(m,r) := \{ t: \; |t| \le r, \; f(z t/\e) \le m \}.
$$
By a crucial lemma of Hal\'asz, the Lebesgue measure of the level sets $|T(m,r)|$
behaves in a regular way (\cite{H}, see \cite{R}, Lemma 3.2):

\begin{lemma}[Regularity]                       \label{Tmr}
  Let $l \in \N$ be such that $l^2 m \le M$. Then
  $$
  |T(m, \frac{\pi}{2})| \le \frac{2}{l} \cdot |T(l^2 m, \pi)|.
  $$
\end{lemma}
Hence, for every $\eta \in (0,1)$ such that $m \le \eta M$, one has:
\begin{equation}                        \label{TmTM}
  |T(m, \frac{\pi}{2})| \le 4 \sqrt{\frac{m}{\eta M}} \cdot |T(\eta M, \pi)|.
\end{equation}
(Apply Lemma~\ref{Tmr} with $l = \lfloor \sqrt{\frac{\eta M}{m}} \rfloor$).

\medskip

Now we can estimate the integral in \eqref{SBP via int} by the
integral distribution formula. Using \eqref{TmTM} for small $m$ and
the trivial bound $|T(m, \pi/2)| \le \pi$ for large $m$, we get
\begin{align}                       \label{SBP via T}
p_\e(a)
  &\le C \sup_{z \ge 1} \int_{-\pi/2}^{\pi/2} \exp \big( - \b f(zt/\e) \big) \; dt \nonumber\\
  &\le C \int_0^\infty |T(m, \frac{\pi}{2})| \ \b e^{- \b m} \; dm \nonumber\\
  &\le C \int_0^{\eta M} 4 \sqrt{\frac{m}{\eta M}} \cdot |T(\eta M, \pi)|
    \ \b e^{- \b m} \; dm
    + C \int_{\eta M}^\infty \pi \ \b e^{- \b m} \; dm \nonumber\\
  &\le \frac{C_1}{\sqrt{\b \eta M}} \cdot |T(\eta M, \pi)|
    + C \pi e^{- \b \eta M} \nonumber\\
  &\le \frac{C_2 B}{\sqrt{\eta n}} \cdot |T(\eta n, \pi)|
    + C \pi e^{- c_2 \eta n / B^2}.
\end{align}
In the last line, we used Lemma~\ref{M} and the definition \eqref{pz} of $\beta$.

\subsection{Recurrence set}

We shall now bound the measure of the level set
$|T(\eta n, \pi)|$ by a quantity of ergodic nature, the density
of the {\em recurrence set} of $a$.

Consider any $t \in T(\eta n, \pi)$ and set $y := z/2\e$.
Then $y \ge 1/2 \e$, and
\begin{equation}                                \label{f as sum}
  f(zt/\e) =  \sum_{k=1}^n \sin^2 (a_k y t) \le \eta n.
\end{equation}
Let us fix
\begin{equation}                                \label{eta}
  \eta := \frac{\a^2 \k}{4n}.
\end{equation}
Then at least $n-\k$ terms in the sum in \eqref{f as sum} satisfy
$$
\sin^2 (a_k y t) \le \frac{\eta n}{\k}
= \frac{\a^2}{4} < \frac{1}{144},
$$
which implies for those terms that
$\dist(a_k y t, \pi \Z) \le \a$.
Thus $yt/\pi$ belongs to the recurrence set of $a$,
which we define as follows:

\begin{definition}[Recurrence set]
  Let $\a \in (0,1)$ and $\k \ge 0$. The {\em recurrence set}
  $I(a) = I_{\a,\k}(a)$ of a vector $a \in \R^n$ is defined as the
  set of all $t \in \R$ such that all except $\k$ coordinates of the vector
  $t a$ are of distance at most $\a$ from $\Z$.
\end{definition}

Regarding $t$ as time, we can think of the recurrence set as the
moments when most of the particles moving along the unit
torus with speeds $a_1,\ldots,a_n$ return close to their initial positions.

\medskip

Our argument thus shows that
$T(\eta n, \pi) \subseteq \frac{\pi}{y}\, I_{\a,\k}(a)$. Thus
$$
|T(\eta n, \pi)|
\le \big| \frac{\pi}{y} \, I_{\a,\k}(a) \cap [-\pi,\pi] \big|
  =  \frac{\pi}{y} \cdot |I_{\a,\k}(a)  \cap [-y,y]|.
$$
The quantity
\[
  \dens(I,y):= \frac{1}{2y} \cdot |I \cap [-y,y]|
\]
can be interpreted as the {\em density} of the set $I$.
We have thus shown that
$$
|T(\eta n, \pi)| \le 2 \pi \, \dens(I_{\a,\k}(a),y).
$$
Using this bound and our choice \eqref{eta} of $\eta$ in \eqref{SBP via T},
we conclude that:
\begin{equation}                                  \label{prob by  I}
p_\e(a)
\le \frac{C_3 B}{\a \sqrt{\k}}
  \cdot \sup_{y \ge 1/2 \e} \dens(I_{\a,\k}(a),y)
  + C \pi e^{- c_3 \a^2 \k / B^2}.
\end{equation}

\subsection{Density of the recurrence set}

It remains to bound the density of the recurrence set $I(a)$
by the reciprocal of the essential LCD $D(a)$. We will derive this from the following
structural lemma, which shows that: (1) the recurrence set has lots of gaps;
(2) each gap bounds below the essential LCD of $a$.

For $t \in \R$, by $[t]$ we denote an integer nearest to $t$.

\begin{lemma}[Gaps in the recurrence set]                   \label{l: scattered}
  Under the assumptions of Theorem~\ref{t: small ball precise},
  let $t_0 \in I_{\a,\k}(a)$.
  Then:
  \begin{enumerate}
    \item $t_0 + 3 \a \not \in I_{\a,\k}(a)$.
    \item Let $t_1 \in I_{\a,\k}(a)$ be such that $t_1 > t_0 + 3 \a$.
    Then $t_1 - t_0 \ge D_{2\a, 2\k}(a)$.
  \end{enumerate}
\end{lemma}
Since $D_{2\a, 2\k}(a) \ge (1 - 2\a)/K > 4 \a$, this lemma implies that
the recurrence set $I$ has gaps of size at least $D_{2\a, 2\k}(a) - 4 \a$.

\begin{proof}
{\em Part 1. }
Since $t_0 \in I_{\a,\k}(a)$, there exists a set
$\s_0 \subseteq \{1 \etc n\}$ of cardinality $|\s_0| \ge n-\k$
and such that for $p_k := [t_0 a_k]$ we have:
\begin{equation}                        \label{sigma 0}
    |t_0 a_k - p_k| \le \a \quad \text{for all } k \in \s_0.
\end{equation}

Let $t := t_0 + 3 \a$. Recall that $1 \le a_k \le K$ for all $k \in \{1 \etc n\}$.
By \eqref{sigma 0}, we have for all $k \in \s_0$:
\begin{align}
t a_k &= t_0 a_k + 3 \a \cdot a_k
  \ge p_k - \a + 3\a
  > p_k + \a;                               \label{sak large}\\
t a_k &\le p_k + \a + 3 \a \cdot a_k
  \le p_k +\a +1/2
  < p_k + 1 -\a.               \notag
\end{align}
In the last inequality, we used the assumption $\a < 1/6K \le 1/6$.
It follows that $\dist(ta_k, \Z) > \a$ for all $k \in \s_0$.
Thus $t \not \in I_{\a,\k}(a)$.
Part 1 is proved.

\medskip

{\em Part 2. }
Since $t_1 \in I_{\a,\k}(a)$, there exists a set
$\s_1 \subseteq \{1 \etc n\}$ of cardinality $|\s_1| \ge n-\k$
and such that for $q_k := [t_1 a_k]$ we have:
\begin{equation}                            \label{sigma 1}
  |t_1 a_k - q_k| \le \a \qquad \text{for all } k \in \s_1.
\end{equation}
Set $\s := \s_0 \cap \s_1$. Then $|\s| \ge n - 2\k$.
Moreover, \eqref{sigma 0} and \eqref{sigma 1} yield:
$$
|(t_1 - t_0)a_k - (q_k - p_k)| \le 2 \a \qquad \text{for all } k \in \s.
$$
Since $t_1 > t$, \eqref{sak large} implies that
\begin{equation}                            \label{s1ak large}
  t_1 a_k > t a_k > p_k + \a \qquad \text{for all } k \in \s.
\end{equation}
Hence, by \eqref{sigma 1} and \eqref{s1ak large}, $q_k - p_k > 0$
for all $k \in \s$.
By the definition of the essential LCD, this means that
\[
  t_1 - t_0 \ge D_{2\a, 2\k}(a).
\]
This completes the proof.
\end{proof}

We can use Lemma~\ref{l: scattered} to bound the density of the
recurrence set via the reciprocal of the essential LCD.

\begin{lemma}[Recurrence set via essential LCD]         \label{rec via LCD}
  Under the assumptions of Theorem~\ref{t: small ball precise},
  we have for every $y > 0$:
  \begin{equation}                          \label{eq rec via LCD}
    \dens(I_{\a,\k}(a),y)
    \le 3 \a \Big( \frac{1}{2y} + \frac{2}{ D_{2\a, 2\k}(a)} \Big).
  \end{equation}
\end{lemma}

\begin{remark}
  The contribution of the first term in \eqref{eq rec via LCD}
  comes from the $O(\a)$-neighborhood of zero, which is contained
  in the recurrence set. This is the initial
  time when all of the moving particles are still close to $0$.
\end{remark}

\begin{proof}
Denote $I := I_{\a,\k}(a) \cap [-y,y]$.
This set is closed and nonempty (it contains $0$). Set
$t_0 := \min \{ t :\; t \in I \}$.
If $I \subseteq [t_0, t_0 + 3 \a]$, then
\begin{equation}                                    \label{I small}
  \dens(I,y) = \frac{|I|}{2y} \le \frac{3 \a}{2y},
\end{equation}
which completes the proof in this case.

Assume then that $I \not \subseteq [t_0, t_0 + 3 \a]$.
Then we can define inductively the maximal sequence
of points $t_1, t_2, \ldots, t_L \in I$ by
$$
t_l := \min \{ t :\; t \in I; \; t > t_{l-1} + 3 \a \}.
$$
Note that by Lemma~\ref{l: scattered}, $t_{l-1} + 3 \a \not\in I$.
Thus the strict inequality in the definition of $t_l$ can be
replaced by the non-strict inequality, so the minimum makes sense.

Part 1 of Lemma \ref{l: scattered} yields
$$
I \subseteq \bigcup_{l=0}^L [t_l, t_l + 3 \a),
$$
while part 2 implies
$$
t_L - t_0 \ge \sum_{l=1}^L (t_l - t_{l-1})
\ge L \cdot D_{2\a, 2\k}(a).
$$
On the other hand,
since $t_0, t_L \in I \subseteq [-y,y]$, we have $t_L - t_0 \le 2y$.
We conclude that
$$
\dens(I,y)
  \le \frac{|\bigcup_{l=0}^L [t_l, t_l + 3 \a)|}{t_L-t_0}
  \le \frac{(L+1)\cdot 3 \a}{L \cdot D_{2\a, 2\k}(a)}
  \le \frac{6 \a}{ D_{2\a, 2\k}(a)}.
$$
This completes the proof.
\endproof

\medskip

By \eqref{prob by I} and Lemma~\ref{rec via LCD}, we conclude that
$$
p_\e(a)
\le \frac{C_4 B}{\sqrt{\k}}
  \Big( \e + \frac{1}{D_{2\a,2\k}(a)} \Big)
  + C \pi e^{- c_3 \a^2 \k / B^2}
$$
for all $\e < \pi/4$ (which was our assumption \eqref{delta small}.)

This completes the proof of Theorem~\ref{t: small ball precise}.

\end{proof}

\subsection{Small ball probability for general coefficients}

In view of the applications, we will state Theorem~\ref{t: small ball intro}
for a general coefficient vector $a$, not necessarily with well comparable
coefficients as in \eqref{K1 K2}. This is easy to do by restricting $a$
onto its spread part, which we define as follows:

\begin{definition}[Spread part]                     \label{d: spread}
  Let $0 < K_1 < K_2$ be fixed.
  For a vector $x \in \R^{n}$, we consider the subset
  $\s(x) \subseteq \{1,\ldots,n\}$ defined as
  $$
  k \in \s(x) \quad \text{if} \quad K_1 \le |n^{1/2} x_k| \le K_2,
  $$
  and, if $\s(x) \ne \emptyset$, we define the {\em spread part of $x$} as
  $$
  \hat{x} := (n^{1/2} x_k)_{k \in \s(x)}.
  $$
  If $\s(x) = \emptyset$, the spread part of $x$ is not defined.
\end{definition}

As an immediate consequence of Restriction Lemma~\ref{l: restriction}
and Theorem~\ref{t: small ball intro}, we obtain:

\begin{corollary}[Small ball probability for general vectors]   \label{SBP general}
  Let $\xi_1,\ldots,\xi_n$ be random variables as in
  Theorem~\ref{t: small ball intro}.
  Let $a \in \R^n$ be a vector of real coefficients whose spread part
  $\hat{a}$ is well defined (for some fixed truncation levels $K_1, K_2 > 0$).
  Let $\a \in (0,1)$ and $\b \in (0,1/2)$.
  Then for for every $\e \ge 0$ one has
  \[
    p_\e(a)
    \le \frac{C}{\sqrt{\b}}
      \Big( \e + \frac{1}{\sqrt{n} \, D_{\a,\b n}(\hat{a})} \Big)
        + C e^{-c \a^2 \b n},
  \]
  where $C, c > 0$ depend (polynomially) only on $B, K_1, K_2$.
\end{corollary}

\begin{remark}
  As a convention throughout the paper, we set $D_{\a,\k}(\hat{a}) = 0$
  if $\hat{a}$ is not defined.
\end{remark}

\begin{remark}
  A small ball probability bound similar to Theorem~\ref{t: small ball precise}
  can be proved with a weaker assumption on the coefficient vector.
  Namely, \eqref{12} can be replaced by
  $$
  \|a\|_1 \ge n, \quad \|a\|_2 \le K \sqrt{n}.
  $$
\end{remark}

\section{Invertibility of random matrices via small ball probability}
\label{s: main proof}

We return here to the invertibility problem for random matrices
that we began to study in Section~\ref{s: median}, and we improve
the Weak Invertibility Theorem~\ref{t: weak} by reducing the
polynomial term $n^{1/2}$ to an exponentially small order $c^n$.

\begin{theorem}[Strong invertibility]                 \label{t: strong}
  Let $\xi_1,\ldots,\xi_n$ be independent centered random variables
  with variances at least $1$ and fourth moments at most $B$.
  Let $A$ be an $n \times n$ matrix whose
  rows are independent copies of the random vector $(\xi_1,\ldots,\xi_n)$.
  Let $K \ge 1$. Then for every $\e \ge 0$ one has
  \begin{equation}              \label{eq: strong}
  \P \big( s_n(A) \le \e n^{-1/2} \big)
    \le C\e + c^n + \P (\|A\| > K n^{1/2}),
  \end{equation}
  where $C> 0$ and $c \in (0,1)$ depend (polynomially) only on $B$ and $K$.
\end{theorem}

This result implies the Subgaussian Invertibility Theorem~\ref{t: subgaussian}:
indeed, the last term in \eqref{eq: strong} is exponentially small
by Lemma~\ref{l: norm subgaussian}.

\medskip

The imprecise term $n^{-1/2}$ in the Weak Invertibility Theorem~\ref{t: weak}
came from from the Weak Distance Bound, Lemma~\ref{l: dist},
which estimated the distance between a random vector and
a random hyperplane.
Thus, in order to complete the proof of the Strong Invertibility
Theorem~\ref{t: strong}, it suffices to improve the bound in
Weak Distance Bound (Lemma~\ref{l: dist}) as follows:

\begin{theorem}[Strong Distance Bound]          \label{t: distance}
  Let $A$ be a random matrix as in Theorem~\ref{t: strong}.
  Let $X_1,\ldots,X_n$ denote its column vectors, and consider
  the subspace $H_n = \Span(X_1,\ldots,X_{n-1})$.
  Let $K \ge 1$. Then for every $\e \ge 0$, one has
  $$
  \P \big( \dist(X_n, H_n) < \e \text{ and } \|A\| \le K n^{1/2}  \big)
  \le C_7(\e + c^n),
  $$
  where $C_7$ and $c \in (0,1)$ depend only on $B$ and $K$.
\end{theorem}

\begin{remark}
  For random vectors with independent $\pm 1$ coordinates, a weaker
  bound $\P \big( \dist(X_n, H_n) < \frac{1}{4n} \big) \le C \log^{-1/2} n$
  was proved by Tao and Vu \cite{TV det}.
\end{remark}

\subsection{Essential LCD of the random normal}

As in Section~\ref{s: dist via SBP}, we shall estimate the
distance by using the the {\em random normal} $X^*$, a unit
normal of the subspace $H_n$. The inequality \eqref{dist}
reduces the problem to a lower bound on $|\< X^*, X_n\> |$.

The random normal $X^*$ is convenient to control via the random
matrix $A'$, the $(n-1) \times n$ matrix with rows
$X_1,\ldots,X_{n-1}$. Thus $A'$ is the submatrix of $A^T$ obtained
by removing the last row. By the definition of the random normal,
$$
A' X^* = 0.
$$
We will use this observation as follows:
\begin{equation}                                \label{normal control}
  \text{If $\|A'x\|_2 > 0$ for all vectors $x$ in some
  set $S$, then $X^* \not\in S$.}
\end{equation}
Thus, a weak (qualitative) invertibility of the random matrix $A'$
on $S$ will help us to ``navigate'' the random normal $X^*$ away from
undesired subsets $S$ of the unit sphere.

We shall use this approach to prove that the essential LCD of the
random normal is exponentially large, with probability exponentially
close to $1$. This will allow us to use the full strength of the Small Ball
Probability Theorem~\ref{t: small ball intro} in order to bound
$|\< X^*, X_n\> |$ from below.

Recall that $\hat{x}$ denotes the spread part of a vector $x$ with
some fixed truncation levels $K_1, K_2$, see Definition~\ref{d:
spread}.

\begin{theorem}[Random normal]                      \label{random normal}
  Let $X_1,\ldots,X_{n-1}$ be random vectors as in Theorem~\ref{t: distance}.
  Consider a unit vector $X^*$ orthogonal to all these vectors.
  Let $K \ge 1$.
  Then there exist constants $K_1, K_2, \a, \b, c, c' > 0$ that depend only
  on $B$ and $K$, and such that
  $$
  \P \big( D_{\a, \b n} (\widehat{X^*}) < e^{cn}
    \text{ and } \|A\| \le K n^{1/2} \big)
  \le e^{-c'n}.
  $$
\end{theorem}

Intuitively, the components of a random vector should be arithmetically
incomparable to the extent that their essential LCD is exponential in $n$.
In the case of the random normal $X^*$, its components are not independent,
and it requires some work to confirm this intuition.

We shall prove that the random matrix $A'$ is likely to be
invertible on the subsets $S_D$ of the unit sphere where the
essential LCD is of order $D$, for each $D$ below an exponential order.
Then, by observation \eqref{normal control}, the random normal $X^*$
will not lie in such $S_D$. Therefore, the essential LCD of $X^*$
will be at least of exponential order.

\subsection{The level sets of the essential LCD}

Fix $K \ge 1$ for the rest of the proof. We shall first
choose the truncation levels $K_1 = K_1(B,K)$, $K_2 = K_2(B,K)$
in the definition of the spread part $\widehat{X^*}$ of the random normal.

Our soft invertibility argument in Section~\ref{s: decomposition}
was based on considering separately compressible and incompressible
vectors, forming the sets
$\Comp = \Comp(\d,\rho)$ and $\Incomp = \Incomp(\d,\rho)$ respectively,
see Definition~\ref{d: compressible}.
The parameters $\d, \rho > 0$ in the definition of these vectors
were chosen in Lemma~\ref{l: compressible}
depending only on $B$ and $K$.

For every incompressible vector $x$, its spread part is
proportionally large. Indeed, by Lemma~\ref{l: spread},
there exist $K_1, K_2, c_0 > 0$ that depend only on $B$ and $K$,
and such that for the truncation levels
$K_1$ and $K_2$ one has $\supp(\hat{x}) \ge c_0 n$.
For the future convenience,
we consider the even integer $n_0 := 2 \lfloor c_0 n / 2 \rfloor$.
Thus we have:
\begin{equation}                                \label{n0}
  \text{Every $x \in \Incomp$ satisfies $|\supp(\hat{x})| \ge n_0
    \ge \frac{c_0}{2}\, n$.}
\end{equation}

We shall choose the value $\a \in (0,1/2)$ later.
By the definition of the essential LCD and of the spread part,
$$
D_{\a, n_0/2}(\hat{x}) \ge (1-\a)/K_2 > 1/2K_2 =: D_0.
$$

\begin{definition}[Level sets of LCD]
  Let $D \ge D_0$. We define the level set $S_D \subseteq S^{n-1}$ as
  $$
  S_D := \{ x \in \Incomp : \; D \le D_{\a, n_0/2}(\hat{x}) < 2D \}.
  $$
\end{definition}

We want to show the invertibility of the random matrix $A'$ on the
level sets $S_D$ for all $D$ up to an exponential order. This will
be done by a covering argument. We will first show the
invertibility on a single vector $x \in S_D$. Next, we will find a
small $(\a/D)$-net in $S_D$. Then, by a union bound, the
invertibility will hold for each point in this net. By
approximation, we will extend the invertibility to the whole $S_D$.

\medskip

The invertibility on a single vector $x \in S_D$ will
easily follow from our general small ball probability
estimates and the Tensorization Lemma~\ref{l: tensorization}.

\begin{lemma}[Invertibility on a single vector]             \label{l: single}
  There exist $c, C_8 > 0$ that depend only on $B$ and $K$,
  and such that the following holds.
  Let $\a \in (0,1)$ and $D_0 \le D < \frac{1}{\sqrt{n}} \, e^{c \a^2 n}$.
  Then for every vector $x \in S_D$ and for every $t \ge 0$, one has
  $$
  \P \big( \|A'x\|_2 < t n^{1/2} \big)
  \le \Big( C_8 t + \frac{C_8}{\sqrt{n} \, D} \Big)^{n-1}.
  $$
\end{lemma}

\begin{proof}
Let $\xi_{k1},\ldots,\xi_{kn}$ denote the $k$-th row of $A'$.
The $k$-th component of $A'x$ is then
$
(A'x)_k = \sum_{j=1}^n x_j \, \xi_{kj} =: \z_k.
$
By Corollary~\ref{SBP general} and by our assumption on $D$,
for every $k$ we have for all $\a \in (0,1)$:
$$
\P( |\z_k| < t)
\le C \Big( t + \frac{1}{\sqrt{n} \, D_{\a,n_0/2} (\hat{x})} \Big)
  + C e^{-c \a^2 n_0 / 2}
\le C' \Big( t + \frac{1}{\sqrt{n} \, D} \Big),
$$
where $C, c, C'$ depend only on $B$ and $K$.

Since $\z_1,\ldots,\z_{n-1}$ are independent random variables
and $\|A'x\|_2^2 = \sum_{k=1}^{n-1} \z_k^2$,
Tensorization Lemma~\ref{l: tensorization} with
$\e_0 =  \frac{1}{\sqrt{n} \, D}$ completes the proof.
\end{proof}

\begin{remark}
  This proof only used the lower bound $D_{\a,n_0/2} (\hat{x}) \ge D$
  in the definition of the level set $S_D$.
\end{remark}

\begin{lemma}[Nets of the level sets]                   \label{l: nets}
  There exist $\a_0 \in (0,1)$, $C_9 > 0$ and $c_9 \in (0,1)$ that depend only on
  $B$ and $K$, and such that the following holds.
  Let $0 < \a < \a_0$ and $D \ge D_0$.
  Then there exists a $(4\a/D)$-net in $S_D$ in the Euclidean metric,
  of cardinality at most
  $$
  \Big( \frac{C_9 D}{\a^{1-c_9}} \Big)^n.
  $$
\end{lemma}

\begin{remark}
  By a simple volumetric estimate (see e.g. \cite{MS}),
  the sphere $S^{n-1}$ has an $\theta$-net of cardinality $(3/\theta)^n$
  for every $\theta > 0$. This implies Lemma~\ref{l: nets} with $c_9 = 0$.
  The fact that the level sets have somewhat smaller cardinality,
  namely with $c_9 > 0$, will be crucial in our argument.
\end{remark}

\begin{proof}
We start by constructing a $(2\a/D)$-net for $S_D$ of the desired
cardinality, whose elements do not necessarily belong to $S_D$.

 Let $x \in S_D$. Recall that
$\supp(\hat{x}) \ge n_0$ by \eqref{n0}. By the definition of
$D(\hat{x}) = D_{\a,n_0/2} (\hat{x})$, there exist $q \in
\R^{\supp(\hat{x})}$ with $n_0/2$ integer coefficients and such that
$$
\|D(\hat{x}) \, \hat{x} - q\|_\infty \le \a.
$$
We can extend $q$ to a vector in $\R^n$ by quantizing its non-integer
coefficients uniformly with with step $\a$. Thus there exists $p \in \R^n$
whose $n_0/2$ coefficients are in $\Z$ and whose other coefficients
are in $\a \Z$, and such that
\begin{equation}                                \label{infty approx}
  \|\sqrt{n} \, D(\hat{x}) \, x - p\|_\infty \le \a.
\end{equation}
(Recall that $\hat{x}$ is a restriction of a vector $\sqrt{n}\, x$).
We thus have $p \in \PP$, where
\begin{equation}                                \label{PP}
  \PP := \bigcup_{|\s| = n_0/2} \Z^\s \oplus \a \Z^{\s^c},
\end{equation}
the union being over all $(n_0/2)$-element subsets $\s$ of $\{1,\ldots,n\}$.

It follows from \eqref{infty approx} and H\"older's inequality that
\begin{equation}                                \label{2 approx}
  \|\sqrt{n} \, D(\hat{x}) \, x - p\|_2 \le \a \sqrt{n}.
\end{equation}
Using $x \in S_D$, we obtain
$$
\Big\| x - \frac{p}{\sqrt{n} \, D(\hat{x})} \Big\|_2
\le \frac{\a}{D(\hat{x})}
\le \frac{\a}{D}
\le \frac{\a_0}{D_0}
\le \frac{1}{4},
$$
if we choose $\a_0 := \min(1, D_0/4)$.

Now we use the following elementary implication, which holds
for every pair of vectors $y$ and $z$ in a Hilbert space: 
if $\|y\| = 1$ and $\|y-z\| \le \d \le 1/4$ then 
$\|y - \frac{z}{\|z\|} \| \le 2 \d$.
This implies
\begin{equation}                                \label{approx}
  \Big\| x - \frac{p}{\|p\|_2} \Big\|_2 \le 2\a/D.
\end{equation}

On the other hand, since $x$ is a unit vector, \eqref{2 approx} implies
$$
\|p\|_2 \le (D(\hat{x}) + \a) \sqrt{n}
  \le 3 \sqrt{n} \, D
$$
where we used that $\a \le \a_0 \le D_0 \le D(\hat{x})$.
We have thus shown that the set
$$
\NN = \Big\{ \frac{p}{\|p\|_2} : \; p \in \PP \cap 3 \sqrt{n} \, D
\cdot B_2^n \Big\} \subset \R^n
$$
is a $(2\a/D)$-net for $S_D$.

Let us estimate the cardinality of $\NN$. There are
$\binom{n}{n_0/2} \le 2^n$ ways to choose the subset $\s$ in
\eqref{PP}. Then
\begin{align*}
|\NN|
  &\le \big| \PP \cap 3 \sqrt{n} \, D \cdot B_2^n \big| \\
  &\le 2^n \cdot \big| \Z^{n_0/2} \cap 3 \sqrt{n} \, D \cdot B_2^{n_0/2} \big|
    \cdot \big| \a \Z^{n-n_0/2} \cap 3 \sqrt{n} \, D \cdot B_2^{n-n_0/2} \big|.
\end{align*}
The Euclidean ball in $\R^d$ of radius $R \sqrt{d}$ and centered at the origin
contains at most $(CR)^d$ integer points, where $C$ is an absolute constant.
Then, using that $n_0 \ge c_0 n/2$, we conclude that
$$
|\NN| \le 2^n \cdot (C \cdot 3D)^{n_0/2} \cdot (C \cdot 3D/\a)^{n-n_0/2}
  \le \Big( \frac{C_9 D}{\a^{1-c_9}} \Big)^n.
$$
Thus, $\NN \subset \R^n$ is a $(2\a/D)$-net for $S_D$ of the
required cardinality. To complete the proof, note that we can make
$\NN$ a subset of $S_D$ using the following standard observation:

\begin{lemma}
 Let $T$ be a metric space and let $E \subset T$.
 Let $\NN \subset T$ be a $\theta$-net of the set $E$.
 Then there exists a $(2\theta)$-net $\NN'$ of $E$
 whose cardinality does not exceed that of $\NN$,
 and such that $\NN' \subset E$.
\end{lemma}
\end{proof}

\begin{remark}
  As we see from \eqref{PP}, we were able to construct a small net
  because of the coarse quantization of a coordinate subspace $\R^\s$
  of proportional dimension, which we could afford due to
  the control of the essential LCD.
  The finer quantization of the complement $\R^{\s^c}$,
  i.e. $\a \Z^{\s^c}$, can be replaced with an arbitrary $\a$-net of that
  subspace. The particular form of the net there does not matter.
\end{remark}

\begin{lemma}[Invertibility on a level set]             \label{l: level}
  There exist $\a, c, c_{10} > 0$ that depend only on $B$ and $K$,
  and such that the following holds.
  Let $D_0 \le D < e^{c n}$.
  Then
  $$
  \P \big( \inf_{x \in S_D} \|A'x\|_2 <  \frac{c_{10}}{D} \, n^{1/2}
    \text{ and } \|A\| \le K n^{1/2} \big)
  \le e^{-n}.
  $$
\end{lemma}

\begin{proof}
Recall that we can assume that $n$ is sufficiently large. We shall therefore
choose a value of $\a$ from the non-empty interval  
$(\frac{1}{\sqrt{n}}, \a_0)$.
Assume that $D_0 \le D < \frac{1}{\sqrt{n}} \, e^{c \a^2 n}$ as 
in Lemma~\ref{l: single}.

We apply Lemma~\ref{l: single} with $t = 5 K \a / D$; thus the
term $C_8 t$ will dominate over the term $C_8/\sqrt{n}\,D$. We
therefore obtain for each $x_0 \in S_D$:
$$
\P \Big( \|A'x_0\|_2 < \frac{5 K \a}{D} \, n^{1/2} \Big) \le
\Big( \frac{C_8'\a}{D} \Big)^{n-1}.
$$
Let $\NN$ be a $(4\a/D)$-net of $S_D$ constructed in Lemma~\ref{l:
nets}. Then taking the union bound, we obtain
\[
  \P \Big( \inf_{x_0 \in \NN} \|A'x_0\|_2 < \frac{5 K \a}{D} \, n^{1/2} \Big)
  \le \Big( \frac{C_9 D}{\a^{1-c_9}} \Big)^n
    \Big( \frac{C_8'\a}{D} \Big)^{n-1}
  \le C_9 D (C_8''\a)^{c_9 n - 1}.
\]
Using the assumption $D < e^{cn}$, we conclude that
\begin{equation}                            \label{on NN}
  \P \Big( \inf_{x_0 \in \NN} \|A'x_0\|_2 < \frac{5 K \a}{D} \, n^{1/2} \Big)
  \le (C_9'\a)^{c_9 n - 1}
  \le e^{-n},
\end{equation}
provided that we choose $\a \ge$ appropriately small 
in the interval $(\frac{1}{\sqrt{n}}, \a_0)$, depending
only on $C_9'$ and $c_9$, which in turn depend only on $B$ and $K$.

We are now ready to bound the event $V$ that $\|A\| \le K n^{1/2}$
and for some $x \in S_D$, $\|A'x\|_2 <  \frac{c_{10}}{D} \, n^{1/2}$.
Assume that $V$ occurs, and choose $x_0 \in \NN$ so that
$\|x-x_0\|_2 \le 4\a/D$. Since $\|A'\| \le \|A\| \le K n^{1/2}$,
we have
$$
\|A'x_0\|_2 \le \|A'x\|_2 + \|A'\| \|x-x_0\|_2
  < \frac{c_{10}}{D} \, n^{1/2} + K n^{1/2} \cdot \frac{4\a}{D}
  \le \frac{5 K \a}{D} \, n^{1/2},
$$
if we choose $c_0 := K \a$ (which thus depends only on $B$ and $K$).
By \eqref{on NN}, this completes the proof.
\end{proof}

\subsection{Proof of the Random Normal Theorem}

Now we prove Theorem~\ref{random normal}.
Let $\a$ and $c$ be as in Lemma~\ref{l: level}.

If $x \in S^{n-1}$ is such that $D(\hat{x}) < e^{cn}$
then, by the definition of the level sets $S_D$,
either $x$ is compressible or $x \in S_D$ for some $D \in \DD$, where
$$
\DD = \{ D :\; D_0/2 \le D < e^{cn}, \; D = 2^k, \; k \in \Z \}.
$$
Therefore, denoting the event that $\|A\| \le K n^{1/2}$ by $U_K$,
we have
$$
\P \big( D(\widehat{X^*}) < e^{cn} \text{ and } U_K \big)
  \le \P \big( X^* \in \Comp \text{ and } U_K \big)
  + \sum_{D \in \DD} \P(X^* \in S_D \text{ and } U_K ).
$$
By Lemma~\ref{normal inc},
$\P( X^* \in \Comp \text{ and } U_K) \le e^{-c_4 n}$.
By \eqref{normal control} and Lemma~\ref{l: level},
for every $D \in \DD$ we have
$$
\P(X^* \in S_D \text{ and } U_K)
  \le \P \big( \inf_{x \in S_D} \|A'x\|_2 =0 \text{ and } U_K \big)
  \le e^{-n}.
$$
Since $|\DD| \le C'n$, we conclude that
$$
\P \big( D(\widehat{X^*}) < e^{cn} \text{ and } U_K \big)
\le e^{-c_4 n} + C'n \cdot e^{-n}
\le e^{-c'n}.
$$
This completes the proof of Theorem~\ref{random normal}. \qed

\subsection{Proof of the Strong Distance Bound and the Strong Invertibility Theorem}

Now we deduce Theorem~\ref{t: distance} from our small ball
probability bound (Corollary~\ref{SBP general}) and the
Random Normal Theorem~\ref{random normal}.

We proceed with a conditioning argument similar to those
used to prove the Weak Distance Bound, Lemma~\ref{l: dist}.
We condition upon a realization of the random vectors
$X_1,\ldots,X_{n-1}$. This fixes realizations of the subspace $H_n$
and the random normal $X^*$. Recall that $X_n$ is independent of
$X^*$. We denote the probability with respect to $X_n$ by $\P_n$,
and the expectation with respect to $X_1,\ldots,X_{n-1}$ by
$\E_{1,\ldots,n-1}$. Then
\begin{multline*}                           \label{strong conditioning}
\P \big( |\< X^*, X_n\> | < \e \text{ and } \|A\| \le K n^{1/2} \big) \\
  \le \E_{1,\ldots,n-1} \P_n \big( |\< X^*, X_n\> | < \e
      \text{ and } D_{\a, \b n}(\widehat{X^*}) \ge e^{cn} \big) \\
    + \P \big( D_{\a, \b n}(\widehat{X^*}) < e^{cn}
      \text{ and } \|A\| \le K n^{1/2}\big).
\end{multline*}

By the Random Normal Theorem~\ref{random normal},
the last term in the right hand side
is bounded by $e^{-c' n}$.
Furthermore, by Corollary~\ref{SBP general}, for any fixed realization
of $X_1,\ldots,X_n$ such that $D_{\a, \b n}(\widehat{X^*}) \ge e^{cn}$
we have
$$
\P_n \big( |\< X^*, X_n\> | < \e \big)
\le C'' \e + C'' e^{-c'' n}.
$$
It follows that
$$
\P \big( |\< X^*, X_n\> | < \e \text{ and } \|A\| \le K n^{1/2} \big)
\le C'' \e + C'' e^{-c'' n} + e^{-c' n}.
$$
By \eqref{dist}, the proof of Theorem~\ref{t: distance} is complete.
\qed

\medskip

Combining Lemma~\ref{l: via distance} and Theorem~\ref{t: distance},
we deduce a strong invertibility bound for a random matrix on
the set of incompressible vectors. This improves a polynomial
term in Lemma~\ref{l: incompressible} to an exponential term:

\begin{lemma}[Strong invertibility for incompressible vectors]
                    \label{l: strong incompressible}
  Let $A$ be a random matrix as in Theorem~\ref{t: strong}.
  Let $K \ge 1$ and $\d,\rho \in (0,1)$.
  Then for every $\e \ge 0$, one has
  $$
  \P \big( \inf_{x \in \Incomp(\d,\rho)} \|A x\|_2 \le \e \rho n^{-1/2} \big)
  \le \frac{C_{11}}{\d} (\e + c^n) + \P (\|A\| > K n^{1/2}),
  $$
  where $C_{11} > 0$ and $c \in (0,1)$ depend only on $B$ and $K$.
\end{lemma}

\medskip

The Strong Invertibility Theorem~\ref{t: strong}
now follows from the decomposition of the sphere \eqref{two terms}
into compressible and incompressible vectors,
and from the invertibility on each of the two parts
established in Lemma~\ref{l: compressible} (see the remark below it)
and Lemma~\ref{l: strong incompressible}
(used for $\d, \rho$ as in Lemma~\ref{l: compressible} and for $\e/\rho$
rather than $\e$). \qed


{\small
\begin{thebibliography}{S 99}

\bibitem{BSY} Z. D. Bai, J. Silverstein, Y. Q. Yin,
  {\em A note on the largest eigenvalue of a large-dimensional sample
  covariance matrix}, J. Multivariate Anal. 26 (1988), 166--168

\bibitem{B} B. Bollob\'{a}s,
  {\em Combinatorics. Set systems, hypergraphs, families of vectors and combinatorial
  probability},
   Cambridge University Press, Cambridge, 1986.

\bibitem{CT} E. J. Candes, T. Tao,
  {\em Near-optimal signal recovery from random projections:
  universal encoding strategies},
  IEEE Trans. Inform. Theory 52 (2004), 5406--5425

\bibitem{DS} K. Davidson, S. J. Szarek,
  {\em Local operator theory, random matrices and Banach spaces},
  Handbook of the geometry of Banach
  spaces, Vol. I,  317--366, North-Holland, Amsterdam, 2001.

\bibitem{Do} D. L. Donoho,
  {\em Compressed sensing},
  IEEE Trans. Inform. Theory 52 (2006), 1289--1306

\bibitem{E 88} A. Edelman,
  {\em Eigenvalues and condition numbers of random matrices},
  SIAM J. Matrix Anal. Appl.  9  (1988), 543--560

\bibitem{E 45} P. Erd\"os,
  {\em On a lemma of Littlewood and Offord},
  Bull. Amer. Math. Soc. 51 (1945), 898--902

\bibitem{E 65} P. Erd\"os,
  {\em Extremal problems in number theory},
  1965 Proc. Sympos. Pure Math., Vol. VIII, pp.181--189
  AMS, Providence, R.I.

\bibitem{Ess} C. G. Esseen,
  {\em On the Kolmogorov-Rogozin inequality for the concentration function},
  Z. Wahrscheinlichkeitstheorie und Verw. Gebiete 5 (1966), 210--216

\bibitem{FF} P. Frankl, Z. Füredi,
   {\em Solution of the Littlewood-Offord problem in high dimensions},
     Ann. of Math. (2)  128  (1988),  no. 2, 259--270.

\bibitem{H 75} G. Hal\'{a}sz,
  {\em On the distribution of additive arithmetic functions},
  Acta Arith. 27 (1975), 143--152

\bibitem{H} G. Hal\'{a}sz,
  {\em Estimates for the concentration function of combinatorial number theory
    and probability},
  Periodica Mathematica Hungarica 8 (1977), 197--211

\bibitem{KKS} J. Kahn, J. Koml\'os, E. Szemer\'edi,
  {\em On the probability that a random $\pm 1$-matrix is singular},
  J. Amer. Math. Soc. 8 (1995), no. 1, 223--240

\bibitem{K 67} J. Koml\'os,
  {\em On the determinant of $(0,\,1)$ matrices},
  Studia Sci. Math. Hungar. 2 (1967), 7--21

\bibitem{Lat} R. Latala,
  {\em Some estimates of norms of random matrices},
  Proc. Amer. Math. Soc. 133 (2005), 1273-1282

\bibitem{LT} M. Ledoux and M. Talagrand,
   {\em Probability in Banach spaces. Isoperimetry and processes},
    Ergebnisse der Mathematik und ihrer Grenzgebiete (3), 23.
    Springer-Verlag, Berlin, 1991.

\bibitem{LPRT}  A.~E.~Litvak, A.~Pajor, M.~Rudelson, N.~Tomczak-Jaegermann,
  {\em Smallest singular value of random matrices and geometry
    of random polytopes},
  Adv. Math. 195 (2005), 491--523

\bibitem{LS} W. V. Li, Q.-M. Shao,
  {\em Gaussian processes: inequalities, small ball probabilities
    and applications}.
  Stochastic processes: theory and methods, 533--597,
  Handbook of Statist., 19, North-Holland, Amsterdam, 2001

\bibitem{MS} V. D. Milman and  G. Schechtman, {\em Asymptotic
    theory of finite-dimensional normed spaces. With an appendix by
    M. Gromov.} Lecture Notes in Mathematics, 1200. Springer-Verlag,
  Berlin, 1986.

\bibitem{O 88} A. M. Odlyzko,
  {\em On subspaces spanned by random selections of $\pm 1$ vectors},
  J. Combin. Theory Ser. A  47  (1988), 124--133

\bibitem{PZ} G. Pan, W. Zhou,
  {\em Circular law, extreme Singular values and potential theory},
  preprint

\bibitem{R} M. Rudelson,
  {\em Invertibility of random matrices: norm of the inverse},
  Annals of Mathematics, to appear

\bibitem{SS} A. S\'ark\"ozy, E. Szem\'eredi,
  {\em \"Uber ein Problem von Erd\"os und Moser},
  Acta Arithmetica 11 (1965), 205--208

\bibitem{S 85} S. Smale,
  {\em On the efficiency of algorithms of analysis},
  Bull. Amer. Math. Soc. (N.S.)  13  (1985), 87--121

\bibitem{So} A. Soshnikov,
  {\em A Note on Universality of the Distribution of the Largest
  Eigenvalues in Certain Sample Covariance Matrices},
  J. Stat. Phys. 108 (2002), 1033--1056

\bibitem{ST} D. Spielman, S.-H. Teng,
  {\em Smoothed analysis of algorithms}.
  Proceedings of the International Congress of Mathematicians,
  Vol. I (Beijing, 2002),  597--606, Higher Ed. Press, Beijing, 2002

\bibitem{Str} D. W. Stroock,
  {\em Probability theory, an analytic view}.
  Cambridge University Press, Cambridge, 1993

\bibitem{Sz}
  S. Szarek,
  {\em Condition numbers of random matrices}, J. Complexity 7 (1991), no. 2,
  131--149.


\bibitem{TV det} T. Tao, V. Vu,
  {\em On random $\pm 1$ matrices: singularity and determinant},
  Random Structures and Algorithms 28 (2006),  1--23

\bibitem{TV singularity} T. Tao, V. Vu,
  {\em On the singularity probability of random Bernoulli matrices},
  J. Amer. Math. Soc., to appear

\bibitem{TV} T. Tao, V. Vu,
  {\em Inverse Littlewood-Offord theorems and the condition
  number of random discrete matrices},
  Annals of Mathematics, to appear

\bibitem{vN} J. von Neumann,
  {\em Collected works.
  Vol. V: Design of computers, theory of automata and numerical analysis}.
  General editor: A. H. Taub. A Pergamon Press Book The Macmillan Co.,
  New York 1963

\bibitem{YBK} Y. Q. Yin, Z. D. Bai, P. R. Krishnaiah,
  {\em On the limit of the largest eigenvalue of the large-dimensional
  sample covariance matrix},
  Probab. Theory Related Fields  78  (1988), 509--521

\end{thebibliography}
}
\end{document}